\documentclass[12pt,twoside,a4paper]{article}
\usepackage{mathrsfs}
\usepackage{mathrsfs}
\usepackage{mathrsfs}
\usepackage{mathrsfs}
\usepackage{mathrsfs}
\usepackage{mathrsfs}
\usepackage{mathrsfs}
\usepackage{stmaryrd}
\usepackage{amsfonts}
\usepackage{amsmath,amssymb}
\usepackage{mathrsfs}
\usepackage{hyperref}
\usepackage{graphicx}
\usepackage{graphicx}

\newtheorem{thm}{Theorem}[section]
\newtheorem{cor}[thm]{Corollary}
\newtheorem{lem}[thm]{Lemma}
\newtheorem{prop}[thm]{Proposition}

\newtheorem{rem}[thm]{Remark}

\numberwithin{equation}{section}\allowdisplaybreaks

\def\leq{\leqslant}
\def\ge{\geqslant}
\def\eps{\varepsilon}
\def\leq{\leqslant}
\def\geq{\geqslant}

\def\Real{{\mathbb{R}}}
\def\NN{{\mathbb{N}}}
\def\FF{{\mathscr{F}^{-1}}}
\def\F{{\mathscr{F}}}
\def\Z{{\mathbb{Z}}}
\newcommand{\norm}[1]{\left\|#1\right\|}
\newcommand{\abs}[1]{\left|#1\right|}

\begin{document}

\pagestyle{myheadings} \markboth{H. Zhang }{GWP for derivative
4NLS}

\title{Global  Well-posedness  for the fourth order nonlinear Schr\"{o}dinger equations with small rough data in high demension}

\author{\bf Hua Zhang\\
{\small \it School of Mathematical Sciences, Peking University,
Beijing 100871, P. R. China}\\
{\small E-mail: zhanghuamaths@163.com(H. Zhang)}}
\date{}
 \maketitle
\thispagestyle{empty}
{\bf Abstract:} For $n\geq 2$, we establish the  smooth effects for
the solutions of the linear fourth order Shr\"{o}dinger equation in
anisotropic Lebesgue spaces with $\Box_k$-decomposition. Using these
estimates, we study the Cauchy problem for the fourth order
nonlinear Schr\"{o}dinger equations with three order derivatives and
obtain the global well posedness for this problem with small data in
modulation space $M^{9/2}_{2,1}({\Real^{n}})$.

 \\

{\bf MSC 2000:} 35 L 30, 46 E 35 47 D 99.

\section{Introduction}
In our earlier paper \cite{Zhang}, we consider the Cauchy problem
for the fourth order nonlinear Schr\"{o}dinger equations with three
order derivatives (4NLS)
\begin{align}
{\rm i} u_{t}+\Delta^{2}u-\eps\Delta u=F((\partial_x^\alpha
u)_{\abs{\alpha}\leq 3},
(\partial_x^\alpha\bar{u})_{\abs{\alpha}\leq 3}),\quad
u(0,x)=u_0(x),
 \label{gNLS}
\end{align}
where $\eps\in\{0,\, 1\}$,  $u$ is a complex valued function of
$(t,x)\in \Real\times \Real^n$.
\begin{align}
\Delta u=-\FF\abs{\xi}^{2}\F u,\quad \Delta^2 u=\FF\abs{\xi}^{4}\F
u,\label{Delta-pm}
\end{align}
$F:\mathbb{C}^{\frac{1}{3}n^3+2n^{2}+\frac{11}{3}n+2}\longrightarrow
\mathbb{C}$ is a polynomial of the form
\begin{align}
F(z)=P(z_{1},...,
z_{\frac{1}{3}n^3+2n^{2}+\frac{11}{3}n+2})=\sum_{m+1\leq\mid\beta\mid\leq
M+1}c_{\beta}z^{\beta},\quad c_{\beta}\in \mathbb{C}.
 \label{poly}
\end{align}
$2+8/n\leq m\leq M, m,M\in\NN$. In this paper, we keep on studying
this problem mainly with the method in \cite{WangNew} .

The fourth order nonlinear Schr\"odinger equation, including its
special forms, arise in  deep water wave dynamics, plasma physics,
optical communications (see \cite{Dys}). A large amount of
 work has been devoted to the  Cauchy problem of
dispersive equations, such as
\cite{Ben,Ca,Chihara,CS,GT,Guo1,Hao,Hao1,Huo,KPV1,KPV2,KPV3,KPV98,KPV4,Seg2003,Seg2004}
and references  therein. In \cite{Seg2003}, by using the method of
Fourier restriction norm, Segata studied a special fourth order
nonlinear Schr\"odinger equation in one dimensional space. And the
results have been improved in \cite{Huo,Seg2004}.

In \cite{WangNew}, Wang, Han and Huang discussed
\begin{align}
{\rm i} u_{t}+\Delta_{\pm}u=F(u,\bar{u},\nabla u,\nabla
 \bar{u}),\quad u(0,x)=u_{0}(x)
 \label{w5}
\end{align}
where $\Delta_{\pm}u=\sum_{i=1}^n \eps_i\partial_{x_i}^2u$ and
$\eps_i\in\{-1,1\}, i=1,...,n$.  They  proved \eqref{w5} is global
well-posed in modulation spaces $M^{s}_{2,1}(\Real^{n}), s\geq3/2$.

\subsection{$M^s_{2,1}$ and $B^s_{2,1}$}

In this paper, we are mainly interested in the cases that the
initial data $u_0$ belongs to the modulation space $ M^s_{2,1}$ for
which the norm can be equivalently defined in the following way (cf.
\cite{Fei2,Wang1,Wang2,Wang3}):
\begin{align}
\|f\|_{M^s_{2,1}}= \sum_{k\in \mathbb{Z}^n} \langle k \rangle^s
\|\mathscr{F} f\|_{L^2(Q_k)},  \label{mod-s21}
\end{align}
where $\langle k\rangle=1+|k|$, $Q_k= \{\xi: -1/2\le \xi_i-k_i <1/2,
\ i=1,...,n\}$. For simplicity, we write $M_{2,1}=M^0_{2,1}$.  Since
only the modulation space $M^s_{2,1}$ will be used in this paper, we
will not state the defination of the general modulation spaces
$M^s_{p,q}$, one can refer to Feichtinger \cite{Fei2}. Modulation
spaces $M^s_{2,1}$ are related to the Besov spaces $B^s_{2,1}$ for
which the norm is defined as follows:
\begin{align}
\|f\|_{B^s_{2,1}}=\|\mathscr{F} f\|_{L^2(B(0,1))}+ \sum^\infty_{j=
1} 2^{sj} \|\mathscr{F} f\|_{L^2(B(0,2^j) \setminus B(0, 2^{j-1}))},
\label{Besov-s21}
\end{align}
where $B(x_0,R):= \{\xi\in \mathbb{R}: \ |\xi-x_0|\le R \}$.  It is
known that there holds the following  optimal inclusions between $
B^{n/2+s}_{2,1}$, $ M^s_{2,1}$ and $B^s_{2,1}$ (cf.
\cite{Toft,SuTo,Wang3}):
\begin{align}
B^{n/2+s}_{2,1}\subset M^s_{2,1} \subset B^s_{2,1}.
\label{Besov-s21}
\end{align}
 So, comparing $M^{s}_{2,1}$ with
$B^{s+n/2}_{2,1}$,  we see that $M^{s}_{2,1}$ contains a class of
functions $u$ satisfying $\|u\|_{M^s_{2,1}}= \infty$ but
$\|u\|_{B^{s+n/2}_{2,1}} \ll 1$. On the other hand, we can also find
a class of rough functions $u$ satisfying $\|u\|_{B^s_{2,1}}=
\infty$ but $\|u\|_{M^{s}_{2,1}} \ll 1$. We have
$B^{n/2}_{2,1}\subset M_{2,1} \subset L^\infty \cap L^2$, this
embedding is also optimal.

\subsection{Main results}

We now give our results, the notations used here can be found in the
section \ref{Notation}.

\begin{thm} \label{DNLS-mod}
Let $n\ge 2$, $ 2\le m \le M <\infty, \; m>8/n$. Assume that $u_0
\in M^{9/2}_{2,1}$ and $\|u_0\|_{M^{9/2}_{2,1}} \le \delta$ for some
small $\delta>0$. Then \eqref{gNLS} has a unique global solution
$u\in C(\mathbb{R}, M^{9/2}_{2,1}) \cap X$, where
\begin{align}
\!\!\!\!\!\! \|u\|_{X} =  &\sum_{\alpha=0,3} \ \sum_{i, \, \ell=1}^n
\ \sum_{k\in \mathbb{Z}^n, \ |k_i|= k_{\rm max} >4} \langle
k_i\rangle^{3} \left\|\partial^\alpha_{x_\ell} \Box_k u
\right\|_{L^\infty_{x_i}
L^2_{(x_j)_{j\not= i}}L^2_t(\mathbb{R}^{1+n})} \nonumber\\
&  + \sum_{\alpha=0,3} \ \sum_{i,\, \ell=1}^n \ \sum_{k\in
\mathbb{Z}^n} \langle k \rangle^{3/2-3/m} \left\|
\partial^\alpha_{x_\ell} \Box_k u \right\|_{L^{m}_{x_i}
L^\infty_{(x_j)_{j\not= i} }
L^\infty_t(\mathbb{R}^{1+n})} \nonumber\\
& + \sum_{\alpha=0,3} \ \sum_{\ell=1}^n \sum_{k\in \mathbb{Z}^n}
\langle k \rangle^{3/2} \left\|
\partial^\alpha_{x_\ell} \Box_k u \right\|_{L^\infty_t L^2_x  \bigcap  L^{2+m}_{x,t}
(\mathbb{R}^{1+n})}, \label{modX.2}
\end{align}
where $k=(k_1,...,k_n)$.  Moreover, $\|u\|_{X} \lesssim \delta$.
\end{thm}
In Theorem \ref{DNLS-mod}, if $u_0 \in M^s_{2,1}$ with $s>9/2$, then
we have $u\in C(\mathbb{R}, M^{s}_{2,1})$.  When the nonlinearity
$F$ has a simple form, say,
\begin{align}
{\rm i} u_t + \Delta^{2}u-\eps\Delta u = \sum^n_{i=1} \lambda_i
\partial^{3}_{x_i} (u^{\kappa_i+1} ), \quad u(0,x)= u_0(x),
\label{DNLS1}
\end{align}
For this special forms, we have

\begin{thm} \label{DNLS1-mod}
Let $n\ge 2$, $\kappa_i \ge  {2\vee \frac{8}{n}}$, $\kappa_i \in
\mathbb{N}$, $\lambda_i \in \mathbb{C}, $ $\kappa=\min_{1\le i \le
n} \, \kappa_i$. Assume that $ u_0 \in M^{3/2}_{2,1}$ and
$\|u_0\|_{M^{3/2}_{2,1}} \le \delta$ for some small $\delta>0$. Then
\eqref{DNLS1} has a unique global solution $u\in C(\mathbb{R},
M^{3/2}_{2,1}) \cap X$, where
\begin{align}
\|u\|_{X_1} =  & \sum^n_{i=1} \sum_{k\in \mathbb{Z}^n, \ |k_i|=
k_{\rm max} >4} \langle k_i\rangle ^{3}  \left\| \Box_k u
\right\|_{L^\infty_{x_i}
L^2_{(x_j)_{j\not=i}}L^2_t(\mathbb{R}^{1+n})} \nonumber\\
&  + \sum^n_{i=1} \sum_{k\in \mathbb{Z}^n} \langle k
\rangle^{3/2-3/\kappa} \left\| \Box_k u \right\|_{L^{\kappa}_{x_i}
L^\infty_{(x_j)_{j\not=i}}L^\infty_t(\mathbb{R}^{1+n})} \nonumber\\
&  +  \sum_{k\in \mathbb{Z}^n}  \langle k \rangle^{3/2} \left\|
\Box_k u \right\|_{L^\infty_t L^2_x  \bigcap  L^{2+\kappa}_{x,t}
(\mathbb{R}^{1+n})}. \label{modX.1}
\end{align}
Moreover, $\|u\|_{X_1} \lesssim \delta$.
\end{thm}
 We remark that in Theorem \ref{DNLS1-mod}, the same result holds
if
 the nonlinear term $
\partial^{3}_{x_i} (u^{\kappa_i+1} )$ is replaced by $
\partial^{3}_{x_i} (|u|^{\kappa_i}u )$ $(\kappa_i\in 2\mathbb{N})$.

\begin{cor} \label{DNLS-cor}
Let $n\ge 2$, $s>(n+3)/2$. Let $X$ and $X_1$ be as in Theorems
\ref{DNLS-mod} and \ref{DNLS1-mod}, respectively. We have the
following results.
 \begin{itemize}
  \item[\rm (i)]  Let $ 2 \le m \le M <\infty$, $m>8/n$. Assume that  $ u_0 \in H^{s+3}$ and $\|u_0\|_{H^{s+1}}
\le \delta$ for some small $\delta>0$. Then \label{4NLS} has a
unique global solution $u\in  X$.
  \item[\rm (ii)]  Let $\kappa_i \ge  {2\vee \frac{8}{n}}$, $\kappa_i \in \mathbb{N}$,
 $\lambda_i \in \mathbb{C}$.
Assume that $ u_0 \in H^{s}$ and $\|u_0\|_{H^{s}} \le \delta$ for
some small $\delta>0$. Then \eqref{DNLS1} has a unique global
solution $u\in X_1$.
\end{itemize}
\end{cor}

\subsection{Notations} \label{Notation}
In this paper, we use the same notation as \cite{WangNew}. The
following are some notations which will be frequently used in this
paper: $\mathbb{C}, \mathbb{R}, \mathbb{N}$ and $ \mathbb{Z}$ will
stand for the sets of complex number, reals, positive integers and
integers, respectively. $c \le 1$, $C>1$ will denote positive
universal constants, which can be different at different places.
$a\lesssim b$ stands for $a\le C b$ for some constant $C>1$, $a\sim
b$ means that $a\lesssim b$ and $b\lesssim a$. We write $a\wedge b
=\min(a,b)$, $a\vee b =\max(a,b)$, $k_{\rm max}=\max_{1\le l \le
n}\abs{k_l}$. We denote by $p'$ the dual number of $p \in
[1,\infty]$, i.e., $1/p+1/p'=1$. We will use Lebesgue spaces
$L^p:=L^p(\mathbb{R}^n)$, $\|\cdot\|_p :=\|\cdot\|_{L^p}$, Sobolev
spaces $H^{s}=(I-\Delta)^{-s/2}L^2$. Some properties of these
function spaces can be found in \cite{Bergh,Triebel}. We will use
the function spaces $L^q_t L^p_x (\mathbb{R}^{n+1})$ and $L^p_x
L^q_t ( \mathbb{R}^{n+1})$ for which the norms are defined by
\begin{align}
& \|f\|_{L^q_{t} L^p_x (\mathbb{R}^{n+1})}= \left\|
\|f\|_{L^{p}_{x}(\mathbb{R}^{n})} \right\|_{L^{q}_t(\mathbb{R})} , \
\ \|f\|_{L^p_x L^q_{t}(\mathbb{R}^{n+1})}= \left\| \|f\|_{
L^{q}_t(\mathbb{R})} \right\|_{L^{p}_{x}(\mathbb{R}^{n})}, \nonumber
\end{align}
$L^p_{x,t}(\mathbb{R}^{n+1}):=L^p_x L^p_{t}(\mathbb{R}^{n+1})$. We
denote by $L^{p_1}_{x_i} L^{p_2}_{(x_j)_{j\not=i}} L^{p_2}_t:=
L^{p_1}_{x_i} L^{p_2}_{(x_j)_{j\not=i}} L^{p_2}_t(\mathbb{R}^{1+n})$
the anisotropic Lebesgue space for which the norm is defined by
\begin{align}
\|f\|_{L^{p_1}_{x_i} L^{p_2}_{(x_j)_{j\not=i}}L^{p_2}_t} =\left\|
\|f\|_{L^{p_2}_{x_1,...,x_{j-1}, x_{j+1},...,x_n}
L^{p_2}_t(\mathbb{R} \times \mathbb{R}^{n-1})}
\right\|_{L^{p_1}_{x_i}(\mathbb{R})} . \label{Notation.1}
\end{align}
It is also convenient to use the notation $L^{p_1}_{x_1}
L^{p_2}_{x_2,...,x_n}L^{p_2}_t: =L^{p_1}_{x_1}
L^{p_2}_{(x_j)_{j\not=1}}L^{p_2}_t $. For any $1<k<n$, we denote by
$\mathscr{F}_{x_1,...,x_k}$ the partial Fourier transform:
\begin{align}
(\mathscr{F}_{x_1,...,x_k} f)(\xi_1,...,\xi_k, x_{k+1},...,x_n) =
\int_{\mathbb{R}^k} e^{-{\rm
i}(x_1\xi_1+...+x_k\xi_k)}f(x)dx_1...dx_k \label{Notation.2}
\end{align}
and by $\mathscr{F}^{-1}_{\xi_1,...,\xi_k}$ the partial inverse
Fourier transform, similarly for $\mathscr{F}_{t,x}$ and
$\mathscr{F}^{-1}_{\tau,\xi}$. $\mathscr{F}:=
\mathscr{F}_{x_1,...,x_n}, \ \mathscr{F}^{-1}:=
\mathscr{F}^{-1}_{\xi_1,...,\xi_n} $. $D^s_{x_i}
=(-\partial^2_{x_i})^{s/2}= \mathscr{F}^{-1}_{\xi_i} |\xi_i|^s
\mathscr{F}_{x_i}$ expresses the partial Riesz potential in the
$x_i$ direction. $\partial^{-1}_{x_i} = \mathscr{F}^{-1}_{\xi_i}
({\rm i} \xi_i)^{-1} \mathscr{F}_{x_i}$. We will use the Bernstein
multiplier estimate; cf. \cite{Bergh,Triebel}. For any $r\in
[1,\infty]$,
\begin{align} \label{Bernstein}
\|\mathscr{F}^{-1}\varphi \mathscr{F} f\|_r\le
C\|\varphi\|_{H^s}\|f\|_r, \quad s>n/2.
\end{align}
We will use the frequency-uniform decomposition operators (cf.
\cite{Wang1,Wang2,Wang3}). Let $\{\sigma_k\}_{k\in \mathbb{Z}^n}$ be
a function sequence satisfying
\begin{align}
 \left\{\begin{array}{l}
\sigma_k(\xi)  \ge c, \quad \forall \; \xi \in Q_k,\\
 {\rm supp}\, \sigma_k \subset \{\xi: |\xi-k| \le \sqrt{n}\},\\
\sum_{k\in \mathbb{Z}^n} \sigma_k (\xi) \equiv 1, \quad \forall \;
\xi \in
\mathbb{R}^n,\\
|D^\alpha \sigma_k(\xi)| \le C_m, \quad \forall \; \xi \in \mathbb{
R}^n,\; |\alpha| \le m \in \mathbb{N}.
\end{array}\right. \label{UD-1}
\end{align}
Denote
\begin{align}
\Upsilon = \left\{ \{\sigma_k\}_{k\in \mathbb{Z}^n}: \;
\{\sigma_k\}_{k\in \mathbb{Z}^n}\;\; \mbox {satisfies \eqref{UD-1}}
\right\}. \label{funct-2}
\end{align}
Let $\{\sigma_k\}_{k\in \mathbb{Z}^n}\in \Upsilon$ be a function
sequence and
\begin{align}
\Box_k := \mathscr{F}^{-1} \sigma_k \mathscr{F}, \quad k\in \mathbb{
Z}^n, \label{funct-3}
\end{align}
which are said to be the frequency-uniform decomposition operators.
One may ask the existence of the frequency-uniform decomposition
operators. Indeed, let $\rho\in \mathscr{S}(\mathbb{R}^n)$ and
$\rho:\, \mathbb{R}^n\to [0,1]$ be a smooth radial bump function
adapted to the ball $B(0, \sqrt{n})$, say $\rho(\xi)=1$ as $|\xi|\le
\sqrt{n}/2$, and $\rho(\xi)=0$ as $|\xi| \ge \sqrt{n} $. Let
$\rho_k$ be a translation of $\rho$: $ \rho_k (\xi) = \rho (\xi- k),
\; k\in \mathbb{Z}^n$.  We write
\begin{align}
\eta_k (\xi)= \rho_k(\xi) \left(\sum_{k\in \mathbb{
Z}^n}\rho_k(\xi)\right)^{-1}, \quad k\in \mathbb{Z}^n.
\label{funct-1}
\end{align}
We have $\{\eta_k\}_{k\in \mathbb{Z}^n}\in \Upsilon$. It is easy to
see that for any $\{\eta_k\}_{k\in \mathbb{Z}^n}\in \Upsilon$,
$$
\|f\|_{M^s_{2,1}} \sim \sum_{k\in \mathbb{Z}^n} \langle k\rangle^s
\|\Box_k f\|_{L^2(\mathbb{R}^n)}.
$$
We will use the function space $\ell^{1,s}_\Box (L^p_{t}L^r_x
(I\times \mathbb{R}^n))$ which contains all of the functions
$f(t,x)$ so that the following norm is finite:
\begin{align} \|f\|_{\ell^{1,s}_\Box (L^p_{t}
L^r_x(I\times \mathbb{R}^n))}:= \sum_{k\in \mathbb{Z}^n} \langle
k\rangle ^s  \|\Box_k f\|_{L^p_{t}L^r_x (I\times \mathbb{R}^n)}.
\label{Mod.7}
\end{align}
For simplicity, we write  $\ell^{1}_\Box (L^p_{t}L^r_x (I\times
\mathbb{R}^n))= \ell^{1,0}_\Box (L^p_{t}L^r_x (I\times
\mathbb{R}^n))$.

This paper is organized as follows. In Section \ref{AGSE} we show
the  smooth effect estimates of the solutions of the  fourth order
linear Schr\"odinger equation in anisotripic Lebesgue spaces with
$\Box_k$-decomposition. In Sections \ref{LE1} and \ref{LE2} we
consider the frequency-uniform localized versions for the global
maximal function estimates, the  smooth effects with
$\Box_k$-decomposition , together with their relations to the
Strichartz estimates. In Sections \ref{pf-thm1} and \ref{pf-thm2} we
prove our Theorems \ref{DNLS1-mod} and \ref{DNLS-mod}, respectively.

\section{Smooth effects with $\Box_k$-decomposition}\label{AGSE}
In this paper,  we always denote
$$
S(t)= e^{{\rm i}t (\Delta^2-\eps\Delta)}= \mathscr{F}^{-1} e^{{\rm
i}t (\abs{\xi}^{4}+\eps\abs{\xi}^{2})}\mathscr{F}, \ \ \mathscr{A} f
(t,x)= \int^t_0 S(t-\tau) f(\tau,x) d\tau.
$$
where  $\eps=0,1 $.
\begin{prop}\label{GSE1}
For  any $k=(k_1,...,k_n)\in
 \mathbb{Z}^n$, $\abs{k_{i}}=k_{\rm max}$, $i=1,...,n$,
 we have
\begin{align}
\norm{\Box_k\partial^{3}_{x_{i}}\mathscr{A}f}_{L^{\infty}_{x_{i}}L^{2}_{(x_{j})j\neq
i}L^{2}_{t}(\Real^{1+n})}\lesssim\norm{\Box_k
f}_{L^{1}_{x_{i}}L^{2}_{(x_{j})j\neq i}L^{2}_{t}(\Real^{1+n})}
\label{smo-eff.1}
\end{align}
{\bf Proof.} Firstly, we assume $\eps=0$. We only give the proof of
the case $i=1$ , the other cases is
identical due to the symmetry. Observing that \\
\begin{align}
\partial^{3}_{x_{1}}\mathscr{A}f=c\mathscr{F}^{-1}_{\tau,\xi}\frac{\xi_{1}^{3}}{\abs{\xi}^{4}-\tau}\F_{t,x}f
\label{smo-eff.2}
\end{align}
Using
 Plancherel's identity, it is equivalent to prove
\begin{align}
\norm{\mathscr{F}^{-1}_{\xi_1}\frac{
\sigma_k(\xi)\xi_{1}^{3}}{\abs{\xi}^{4}-\tau}\mathscr{F}_{x_{1}}f}_{L^{\infty}_{x_{1}}L^{2}_{(\xi_{j})j\neq
1}L^{2}_{t}(\Real^{1+n})}\lesssim \norm{
\mathscr{F}^{-1}_{\xi_1}\sigma_k(\xi)
f}_{L^{1}_{x_{1}}L^{2}_{(\xi_{j})j\neq 1}L^{2}_{\tau}(\Real^{1+n})}
\label{smo-eff.3}
\end{align}
Using Young's inequality, it is suffices to prove
\begin{align}\label{smo-eff.4}
\sup_{x_{1}, \tau,
\xi_{j}(j\neq1)}\abs{\mathscr{F}^{-1}_{\xi_1}\frac{\sigma_k(\xi)\xi_{1}^{3}}{\abs{\xi}^{4}-\tau}}\lesssim
C.
\end{align}
We give the proof of \eqref{smo-eff.4} according to $\tau > 0$ or
$\tau \leq 0$. Obsering that in this case we have $\abs{\xi_{1}}\sim
\max\abs{\xi_{j}},j=1,...,n,
 \abs{\abs{\xi}-k}\leq \sqrt{n}$. Therefore, when $\tau \leq 0$ we
have
\begin{align*}
&&\sup_{x_{1}, \tau,
\xi_{j}(j\neq1)}\abs{\mathscr{F}^{-1}_{\xi_1}\frac{\sigma_k(\xi)\xi_{1}^{3}}{\abs{\xi}^{4}-\tau}}
\lesssim\abs{\int_{\abs{\xi_{1}}\sim k}\frac{1}{\xi_{1}}d\xi_{1}}
\lesssim C
\end{align*}
When $\tau>0$, for simplicity we drop $\sigma_k(\xi)$. From
$\abs{\xi}^{4}-\tau=(\abs{\xi}^{2}-\sqrt{\tau})(\abs{\xi}^{2}+\sqrt{\tau})$,
and let
$\abs{\overline{\xi}}^{2}=\sum\nolimits^{n}_{j=2}{\xi_{j}}^{2}$. We
have
\begin{eqnarray*}
&&\int\frac{\xi_{1}^{3}}{\abs{\xi}^{4}-\tau}e^{ix_{1}\xi_{1}}d\xi_{1}\\
&=&\int\frac{\xi_{1}^{3}}{(\abs{\xi}^{2}-\sqrt{\tau})(\abs{\xi}^{2}+\sqrt{\tau})}e^{ix_{1}\xi_{1}}d\xi_{1}\\
&=&\int\frac{\xi_{1}^{3}}{(\xi_{1}^{2}+\abs{\overline{\xi}}^{2}-\sqrt{\tau})(\xi_{1}^{2}+\abs{\overline{\xi}}^{2}+\sqrt{\tau})}e^{ix_{1}\xi_{1}}d\xi_{1}
\end{eqnarray*}
When $\abs{\overline{\xi}}^{2}-\sqrt{\tau}\geq0$, we easily get the
desired result.\\
When $\abs{\overline{\xi}}^{2}-\sqrt{\tau}<0$, we let
$A^{2}=-(\abs{\overline{\xi}}^{2}-\sqrt{\tau})$ and
$B^{2}=\abs{\overline{\xi}}^{2}+\sqrt{\tau}$. We get
\begin{eqnarray}\label{eps=0}
&&\int\frac{\xi_{1}^{3}}{(\xi_{1}^{2}+\abs{\overline{\xi}}^{2}-\sqrt{\tau})(\xi_{1}^{2}+\abs{\overline{\xi}}^{2}+\sqrt{\tau})}e^{ix_{1}\xi_{1}}d\xi_{1}\\
&=&\int\frac{\xi_{1}}{\xi_{1}^{2}-A^{2}}\frac{\xi_{1}^{2}}{\xi_{1}^{2}+B^{2}}e^{ix_{1}\xi_{1}}d\xi_{1}\nonumber\\
&=&\frac{1}{2}\int(\frac{1}{\xi_{1}+A}+\frac{1}{\xi_{1}-A})\frac{\xi_{1}^{2}}{\xi_{1}^{2}+B^{2}}e^{ix_{1}\xi_{1}}d\xi_{1}\nonumber\\
&=& I+II\nonumber
\end{eqnarray}
For part I, we have
\begin{eqnarray*}
I&=&\frac{1}{2}\int\frac{1}{\xi_{1}+A}\frac{\xi_{1}^{2}+B^{2}-B^{2}}{\xi_{1}^{2}+B^{2}}e^{ix_{1}\xi_{1}}d\xi_{1}\\
&=&\frac{1}{2}\int\frac{1}{\xi_{1}+A}e^{ix_{1}\xi_{1}}d\xi_{1}+\frac{1}{2}\int\frac{1}{\xi_{1}+A}\frac{-B^{2}}{\xi_{1}^{2}+B^{2}}e^{ix_{1}\xi_{1}}d\xi_{1}\\
&=&I_{1}+I_{2}
\end{eqnarray*}
The part $I_{1}$ is bounded onwing to Hilbert transform. For
$I_{2}$, by changes of variables, it suffices to show
\begin{align}
\sup_{x_1}\int\frac{1}{1+\xi_{1}}\frac{1}{1+D^{2}\xi_{1}^{2}}e^{ix_{1}\xi_{1}}d\xi_{1}\lesssim
C
\end{align} where $D=\frac{A}{B}$.
Using the fact that $\F(e^{-\abs{x}})(\xi)=C
\frac{1}{1+\abs{\xi}^{2}}$, we have
\begin{eqnarray*}
\norm{\int\frac{1}{1+\xi_{1}}\frac{1}{1+D^{2}\xi^{2}_{1}}e^{ix_{1}\xi_{1}}d\xi_{1}}_{L^{\infty}_{x_{1}}}&=&\norm{\mathscr{F}^{-1}_{\xi_1}(\frac{1}{1+\xi_{1}})\ast\mathscr{F}^{-1}_{\xi_1}(\frac{1}{1+D^{2}\xi^{2}_{1}})}_{L^{\infty}_{x_{1}}}\\
&\lesssim&\norm{\mathscr{F}^{-1}_{\xi_1}(\frac{1}{1+\xi_{1}})}_{L^{\infty}_{x_{1}}}\norm{\mathscr{F}^{-1}_{\xi_1}(\frac{1}{1+D^{2}\xi^{2}_{1}})}_{L^{1}_{x_{1}}}\\
&\lesssim&\norm{\mathscr{F}^{-1}_{\xi_1}(\frac{1}{1+\xi_{1}})}_{L^{\infty}_{x_{1}}}\norm{\frac{1}{D}e^{-\frac{x_{1}}{D}}}_{L^{1}_{x_{1}}}\\
&\lesssim& C
\end{eqnarray*}
The part $II$ is similar to $I$, so we get the result desired. Now
we consider the case $\eps=1$. Comparing the  proof of the case
$\eps=0$, it suffices to show

\begin{align*}
\sup_{x_{1}, \tau,
\xi_{j}(j\neq1)}\abs{\mathscr{F}^{-1}_{\xi_1}\frac{\sigma_k(\xi)\xi_{1}^{3}}{\abs{\xi}^{4}+\abs{\xi}^{2}-\tau}}\lesssim
C.
\end{align*}
When $\tau \leq 0$, the proof is identical to the case $\eps=0$.
Observing that when $\tau >0$, we can choose
$\tau_{2}=-\frac{1}{2}+\sqrt{\frac{1}{4}+\tau} > 0$ such that
\begin{align}
\abs{\xi}^{4}+\abs{\xi}^{2}-\tau=(\abs{\xi}^{2}-\tau_{2})(\abs{\xi}^{2}+\tau_{2}+1)
\end{align}
which is turn to  \eqref{eps=0}. $\hfill\Box$
\end{prop}

\begin{rem} \rm We assume $|k_{i}|= k_{\rm max}$ in Prop \ref{GSE1}. For the general case, see
 Section \ref{LE2} for details.
\end{rem}
\begin{prop}\label{GSE2}
For any $k=(k_1,...,k_n)\in
 \mathbb{Z}^n$, $\abs{k_{i}}=k_{\rm max}$, $i=1,...,n$, we have
\begin{align}\label{smo-eff.21}
\norm{\Box_kD^{3/2}_{x_{i}}S(t)u_{0}}_{L^{\infty}_{x_{i}}L^{2}_{(x_{j})j\neq
i}L^{2}_{t}(\Real^{1+n})}\lesssim\norm{\Box_k u_{0}}_{L^{2}}
\end{align}
{\bf Proof.} As Prop \ref{GSE1}, we only need to prove the case
$i=1$. By Plancherel's identity, we have
\begin{align*}
& \norm{\Box_k
D^{3/2}_{x_{1}}S(t)u_{0}}_{L^{\infty}_{x_{1}}L^{2}_{(x_{j})j\neq
1}L^{2}_{t}(\Real^{1+n})}\\
& =\norm{
\int\sigma_k(\xi)\abs{\xi_{1}}^{\frac{3}{2}}e^{\rm
it(\abs{\xi}^{4}+\eps\abs{\xi}^{2})}\widehat{u_{0}}(\xi)e^{\rm
ix_{1}\xi_{1}}d\xi_{1}}_{L^{\infty}_{x_{1}}L^{2}_{(\xi_{j})j\neq
1}L^{2}_{t}(\Real^{1+n})}
\end{align*}
We can assume $\xi_{1}>0$, otherwise we let $\xi_{1}'=-\xi_{1}$.
Making variables change $\eta = \abs{\xi}^{4} + \eps\abs{\xi}^{2}$
and using Plancherel's identity, we have
\begin{eqnarray*}
&&\norm{\int\sigma_k(\xi)\xi_{1}^{\frac{3}{2}}e^{\rm i
t\abs{\xi}^{4}}\widehat{u_{0}}(\xi)e^{\rm i
x_{1}\xi_{1}}d\xi_{1}}_{L^{\infty}_{x_{1}}L^{2}_{(\xi_{j})j\neq
1}L^{2}_{t}(\Real^{1+n})}\\
&\lesssim&\norm{\int\sigma_k(\xi(\eta))\xi_{1}^{\frac{3}{2}}(\eta)
e^{\rm i t\eta}\widehat{u_{0}}(\xi(\eta))e^{\rm i
x_{1}\xi_{1}(\eta)}(\abs{\xi}^2+\eps)^{-1}\xi_{1}^{-1}(\eta)d\eta}_{L^{\infty}_{x_{1}}L^{2}_{t}L^{2}_{(\xi_{j})j\neq
1}(\Real^{1+n})}\\
&\lesssim &
\norm{\sigma_k(\xi(\eta))\xi_{1}^{\frac{1}{2}}(\eta)\widehat{u_{0}}(\xi(\eta))(\abs{\xi}^2+\eps)^{-1}}_{L^{2}_{\eta}L^{2}_{(\xi_{j})j\neq
1}(\Real^n)}\\
&\lesssim & \norm{\sigma_k(\xi)\xi_{1}^{\frac{1}{2}}
\widehat{u_{0}}(\xi)(\abs{\xi}^2+\eps)^{-1}(\abs{\xi}^2+\eps)^{1/2}\xi_{1}^{\frac{1}{2}}}_{L^{2}_{\xi_{1}}L^{2}_{(\xi_{j})j\neq
1}(\Real^n)}\\
&=&\norm{\sigma_k(\xi)\xi_{1}(\abs{\xi}^2+\eps)^{-1/2}
\widehat{u_{0}}(\xi)}_{L^{2}_{\xi_{1}}L^{2}_{(\xi_{j})j\neq
1}(\Real^n)}\\
&\lesssim&\norm{u_{0}}_{L^2}
\end{eqnarray*}
$\hfill\Box$
\end{prop}

By the duality of \eqref{smo-eff.21}, we have the following
\begin{prop}\label{GSE3}
For  any $k=(k_1,...,k_n)\in
 \mathbb{Z}^n$, $\abs{k_{i}}= k_{\rm max}$, $i=1,...,n$, we have
\begin{align}
\norm{\Box_k
\partial^{3}_{x_{i}}\mathscr{A}f}_{L^{\infty}_{t}L^{2}_{x}(\Real^{1+n})}\lesssim\norm{ \Box_k D^{3/2}_{x_{i}}f}_{L^{1}_{x_{i}}L^{2}_{(x_{j})j\neq
i}L^{2}_{t}(\Real^{1+n})}.
\end{align}
\end{prop}

\section{Linear estimates with $\Box_k$-decomposition} \label{LE1}

In this section we consider the smooth effect estimates, the maximal
function estimates, the Strichartz estimates and their interaction
estimates for the solutions of the  fourth order Schr\"odinger
equations by using the frequency-uniform decomposition operators.
For convenience, we will use the following function sequence
$\{\sigma_k\}_{k\in \mathbb{Z}^n}$. Now we recall some results in
\cite{WangNew}.

\begin{lem}\label{eq-mod-funct}
Let $\eta_k : \mathbb{R} \to [0,1]$ $(k\in \mathbb{Z})$ be a
smooth-function sequence satisfying condition \eqref{UD-1}. Denote
\begin{align}
\sigma_k (\xi):= \eta_{k_1}(\xi_1)...\eta_{k_n}(\xi_n), \ \
k=(k_1,...,k_n). \label{UD-2}
\end{align}
Then we have $\{\sigma_k\}_{k\in \mathbb{Z}^n} \in \Upsilon$.
\end{lem}

\begin{lem}\label{discret-deriv}
For any $\sigma \in  \mathbb{R}$ and  $k=(k_1,...,k_n)\in
\mathbb{Z}^n$ with $|k_i| \ge 4$, we have
\begin{align}
\|\Box_k D^\sigma_{x_i} u\|_{L^{p_1}_{x_1} L^{p_2}_{x_2,...,x_n}
L^{p_2}_t (\mathbb{R}^{1+n}) } & \lesssim \langle k_i
\rangle^{\sigma} \|\Box_k u\|_{L^{p_1}_{x_1} L^{p_2}_{x_2,...,x_n}
L^{p_2}_t (\mathbb{R}^{1+n}) }. \label{equivalent}
\end{align}
Replacing $D^\sigma_{x_i}$ by $\partial^\sigma_{x_i}$ ($\sigma\in
\mathbb{N}$), the above inequality holds for all $k\in
\mathbb{Z}^n$.
\end{lem}

\begin{rem}\label{eqrm} \rm From the proof, we can get that the $\lesssim$  in Lemma \ref{discret-deriv}
can be strengthened as $\sim$. In fact, taking $v=D^\sigma_{x_i} u$
in \eqref{equivalent} and noticing $\sigma \in \mathbb{R}$, we can
get the reversal inequality of \eqref{equivalent}.
\end{rem}

The next lemma is essentially known, see \cite{Triebel,Wang1}.

\begin{lem}\label{qnorm:pnorm}
Let $\Omega\subset {\Bbb R}^n$ be a compact set with ${\rm
diam}\,\Omega <2R$, $0<p\le q \le \infty.$ Then there exists a
constant $C>0$, which depends only on $p,q$ such that
\begin{align}
\| f\|_q \le C R^{n(1/p-1/q)}\|f\|_p, \quad \forall\; f\in
L^p_\Omega, \nonumber
\end{align}
where $L^p_\Omega =\{f\in {\cal S}'({\Bbb R}^n): {\rm supp}
\hat{f}\subset \Omega, \; \|f\|_p<\infty\}.$
\end{lem}

Here we emphasis that we can find a  constant $C>0$ uniformly holds
for all $k\in {\Bbb Z}^n$ in Lemma \ref{qnorm:pnorm}.

It is known that $S(t)$ satisfy the following $L^{p} - L^{p'}$
estimate:
\begin{align}
\norm{S(t)f}_{p}\lesssim \abs{t}^{-n/4(1-2/p)}\norm{f}_{p'}, \abs{t}
\geq 1 ; 2\leq p\leq\infty,  \label{pp'estimate1}
\end{align}
Using the same procedure as  in \cite{Wang1}, we have
\begin{align}
\norm{\Box_k S(t)f}_{p}\lesssim \sum_{l\in \Lambda}\norm{\Box_{k+l}
f}_{p'}, 2\leq p \leq\infty \label{pp'estimate2}
\end{align}
where $\Lambda=\{l\in \Z^{n}: B(0,\sqrt{n})\cap B(l,\sqrt{n})\neq
\phi\}$. \\
Combining \eqref{pp'estimate1} and \eqref{pp'estimate2}, we have
\begin{align}
\norm{\Box_k S(t)f}_{p}\lesssim (1+\abs{t})^{-n/4(1-2/p)}\sum_{l\in
\Lambda}\norm{\Box_{k+l} f}_{p'}, 2\leq p \leq\infty
\label{pp'estimate3}
\end{align}
Using \eqref{pp'estimate3} and following the procedure in
\cite{Wang2}, we  get the following
\begin{lem}\label{Strichartz-mod}
 Let $2\le p< \infty$, $\gamma \ge 2\vee \gamma(p)$,
\begin{align}
 \frac{4}{\gamma(p)}=n\Big(\frac{1}{2}-\frac{1}{p}\Big). \nonumber
\end{align}
 Then we have
\begin{align}
\left\|S(t) \varphi \right\|_{\ell^1_\Box (L^\gamma(\mathbb{R},
L^p(\mathbb{R}^n))
)} & \lesssim \| \varphi\|_{M_{2, 1}(\mathbb{R}^n)},  \nonumber\\
\left\|\mathscr{A} f \right\| _{\ell^1_\Box (L^\gamma (\mathbb{R},
L^p(\mathbb{R}^n))) \cap \ell^1_\Box (L^\infty(\mathbb{R},
L^2(\mathbb{R}^n))) } & \lesssim \|f\|_{\ell^1_\Box
(L^{\gamma'}(\mathbb{R}, L^{p'}(\mathbb{R}^n)))}. \nonumber
\end{align}
In particular, if $2+8/n\le p< \infty$, then we have
\begin{align}
\left\|S(t) \varphi \right\|_{\ell^1_\Box
(L^p_{t,x}(\mathbb{R}^{1+n})
)} & \lesssim \| \varphi\|_{M_{2, 1}(\mathbb{R}^n)},  \nonumber\\
\left\|\mathscr{A} f \right\| _{\ell^1_\Box
(L^p_{t,x}(\mathbb{R}^{1+n})) \, \cap \, \ell^1_\Box
(L^\infty_tL^2_x(\mathbb{R}^{1+n})) } & \lesssim \|f\|_{\ell^1_\Box
(L^{p'}_{t,x}(\mathbb{R}^{1+n}))}. \nonumber
\end{align}
\end{lem}

In \cite{Io-Ke}, when $n\ge 3$  Ionescu and Kenig  showed the
following maximal function estimates:
\begin{align}
 \|\triangle_k S(t)
u_0\|_{L^{2}_{x_i} L^\infty_{(x_j)_{j\not=i}} L^\infty_t
(\mathbb{R}^{1+n}) } & \lesssim  2^{(n-1)k/2} \| \triangle_k
u_0\|_{L^2(\mathbb{R}^n)}. \label{MFM-2a}
\end{align}
Combining their idea and the frequency-uniform decomposition
operators as \cite{WangNew}, we  obtain the following

\begin{prop}\label{MaxFunct-Mod}
 Let $8/n< q \le \infty$, $q\ge 2$ and  $k=(k_1,...,k_n)\in
 \mathbb{Z}^n$,  we have
\begin{align}
 \|\Box_k S(t)
u_0\|_{L^{q}_{x_i} L^\infty_{(x_j)_{j\not=i}} L^\infty_t
(\mathbb{R}^{1+n}) } & \lesssim  \langle k_{\rm max} \rangle^{3/q}
\| \Box_k u_0\|_{L^2(\mathbb{R}^n)}. \label{MFM-2}
\end{align}
\end{prop}

{\bf Proof.} It suffices to prove
\begin{align}
\norm{\int_{\Real^{n}}\delta_{k}(\xi)e^{{\rm i}x \xi}e^{{\rm
 i}t(|\xi|^4 + \eps|\xi|^2)}u_{0}d\xi}_{L^{q}_{x_1} L^\infty_{(x_j)_{j\not=1}} L^\infty_t
(\mathbb{R}^{1+n})} \lesssim \langle k_{\rm max} \rangle^{3/q} \|
\Box_k u_0\|_{L^2(\mathbb{R}^n)}
 \label{kw1}
 \end{align}
By a standard $TT^*$ method, it suffices to prove
\begin{align}
\norm{\int_{\Real^{n-1}\times \Real}\delta_{k}(\xi)e^{{\rm i}x_{1}
\xi_{1}}e^{{\rm i}\overline{x} \overline{\xi}}e^{{\rm
 i}t(|\xi|^4 + \eps|\xi|^2)}u_{0}(\xi)d\xi_{1}d\overline{\xi}}_{L^{q/2}_{x_1}
 L^\infty_{\overline{x}, t}
(\mathbb{R}^{1+n})} \lesssim \langle k_{\rm max} \rangle^{6/q}
 \label{kw2}
 \end{align}
 where $\overline{k}=(k_{2},... ,k_{n}), \overline{x}=(x_{2},... ,x_{n}), \overline{\xi}=(\xi_{2},...
 ,\xi_{n})$.
 For convenience, we give the details of \eqref{kw2} when $n=2$.
 The general case can be treated similarly.

 When $n=2$, we can write \eqref{kw2} as following
\begin{align}
\norm{\int_{\Real \times \Real}\delta_{k}(\xi)e^{{\rm i}x_{1}
\xi_{1}}e^{{\rm i}x_{2} \xi_{2}}e^{{\rm
 i}t(|\xi|^4 + \eps|\xi|^2)}u_{0}(\xi)d\xi_{1}d\xi_{2}}_{L^{q/2}_{x_1}
 L^\infty_{x_{2}, t}
(\mathbb{R}^{2})} \lesssim \langle k_{\rm max} \rangle^{6/q}
 \label{kw3}
 \end{align}

Making variable changes $\xi_{j}-k_{j}=\mu_{j}, j=1,2$, ie. $\xi -k
=\mu$, we have
\begin{eqnarray}
&&\norm{\int_{\Real \times
\Real}\delta_{0}(\mu)\prod^{2}_{j=1}e^{{\rm i}x_{j} k_{j}}e^{{\rm
i}x_{j} \mu_{j}}e^{{\rm
 i}t(|\xi+k|^{4}+\eps|\xi+k|^{2})}u_{0}(\mu+k)d\mu_{1}d\mu_{2}}_{L^{q/2}_{x_1}
 L^\infty_{x_{2}, t}
(\mathbb{R}^{2})}\nonumber\\
 &&\lesssim \langle k_{\rm max} \rangle^{6/q}
 \label{kw4}
 \end{eqnarray}

 We use the $L$ denotes the left of \eqref{kw4} in later proof.
Expanding the term $|\xi+k|^{4}+\eps|\xi+k|^{2}$, we obtain three
sorts of terms.

Firstly, the term such as  $k_{1}^{4}, k_{2}^{4},
...,2(k_{1}k_{2})^{2}$ have no relation with the integral variable
$\mu$. For these terms, we have $\abs{e^{{\rm i}t
k_{1}^{4}}}=\abs{e^{{\rm i}t k_{2}^{4}}}=...=\abs{e^{2{\rm i}t
(k_{1}k_{2})^{2}}}=1$.

Secondly, noticing $\mu\sim 0$, we can treat terms such as
$\mu_{1}^{4}, \mu_{2}^{4}, (\mu_{1}\mu_{2})^{2}$ as the first case.

Finally,  the main contribution  to the $L$ in the remainder terms
 such as $\mu_{1}k_{1}^{3}, \mu_{2}k_{2}^{3}$.
We treat this case according to $\abs{k_{1}} = k_{\rm max}$ or
$\abs{k_{2}} = k_{\rm max}$. When $\abs{k_{1}} = k_{\rm max}$, for
\eqref{kw4}, we need to prove
\begin{align}
L\lesssim \langle k_{1} \rangle^{6/q}
\label{kw5}
\end{align}

Using the above analysis about $|\xi+k|^{4}+\eps|\xi+k|^{2}$, we
have
\begin{align}
L\approx \norm{\int_{\Real \times
\Real}\delta_{0}(\mu)\prod^{2}_{j=1}e^{{\rm i}x_{j} k_{j}}e^{{\rm
i}x_{j} \mu_{j}}e^{{\rm
 i}t(\mu_{1}k_{1}^{3}+\mu_{2}k_{1}^{3})}u_{0}(\mu+k)d\mu_{1}d\mu_{2}}_{L^{q/2}_{x_1}
 L^\infty_{x_{2}, t}
(\mathbb{R}^{2})} \label{kw6}
\end{align}
In view of Lemma \ref{eq-mod-funct}, we can write $\Box_k=
\mathscr{F}^{-1} \eta_{k_1}(\xi_1)...\eta_{k_n}(\xi_n)\mathscr{F}:=
\mathscr{F}^{-1}
\eta_{k_1}(\xi_1)\eta_{\bar{k}}(\bar{\xi})\mathscr{F}$. For one
thing, in view of  the decay of  $\Box_k S(t)$
\begin{align}
& \|\mathscr{F}^{-1}_{\bar{\xi}} e^{{\rm i} t(|\bar{\xi}|^4 + \eps
|\bar{\xi}|^2 )} \eta_{\bar{k}}(\bar{\xi}) \|_{L^\infty_{\bar{x}}
(\mathbb{R}^{n-1})}
\lesssim (1+|t|)^{-(n-1)/4}, \nonumber\\
& \|\mathscr{F}^{-1}_{\xi_1} e^{{\rm i} t (\xi_1^4 +\eps \xi_1^2)}
\eta_{k_1}(\xi_1) \|_{L^\infty_{x_1} (\mathbb{R})} \lesssim
(1+|t|)^{-1/4}. \nonumber
\end{align}

 For other, integrating by part when $x_{1}\geq4
k_{1}^{3}\abs{t}+1$,
\begin{align}
\abs{\mathscr{F}^{-1}_{\xi_1} e^{{\rm i} t (\xi_1^4 +\eps \xi_1^2)}
\eta_{k_1}(\xi_1)}\lesssim \abs{x_{1}}^{-2}. \nonumber
\end{align}

Therefore, when $x_{1}\geq4 k_{1}^{3}\abs{t}+1$  we have
\begin{align}
L\lesssim\norm{\abs{x_{1}}^{-2}}_{L^{q/2}_{x_{1}}(\Real)}\lesssim C
\end{align}

When $x_{1}\leq4 k_{1}^{3}\abs{t}+1$, observing
$(1+\abs{t})^{-n/4}\lesssim \langle k_1^{3}\rangle^{n/4}(\langle
k_1^{3}\rangle+\abs{x_{1}})^{-n/4}$, we have
\begin{eqnarray}
L&&\lesssim\langle k_1^3\rangle^{1/2}\norm{(\langle
k_1\rangle^3+\abs{x_{1}})^{-n/4}}_{L^{q/2}_{x_{1}}(\Real)}\nonumber\\
&&\lesssim\langle k_1\rangle^{6/q} \label{kw7}
\end{eqnarray}
Observing the above argument also holds with change the place of
$k_{1}$ and $k_{2}$ when $|k_{2}|= k_{\rm max}$ . From this, we can
see why the right of \eqref{MFM-2} is $k_{\rm max}$. $\hfill \Box$

By duality of Proposition \ref{MaxFunct-Mod}, we have  the following
\begin{prop}\label{MaxFunct-Mod-Dual}
For $2\le  q \le \infty$, $q>8/n$ and $k=(k_1,...,k_n)\in
 \mathbb{Z}^n$, we have
\begin{align}
\left \|\Box_k \int_\mathbb{R} S(t-\tau) f(\tau) d\tau
\right\|_{L^\infty_t L^2_x(\mathbb{R}^{1+n})} & \lesssim \langle
k_{\rm max} \rangle^{3/q} \| \Box_k f\|_{L^{q'}_{x_i}
L^1_{(x_j)_{j\not=i}} L^1_t (\mathbb{R}^{1+n}) }. \label{MFM-2-dual}
\end{align}
\end{prop}

In view of Propositions \ref{GSE1} and \ref{GSE3}, we have
\begin{prop}\label{1order-sm-mod}
For any $k=(k_1,...,k_n)\in
 \mathbb{Z}^n$ and $\abs{k_{i}} = k _{max}, i=1,...,n$, we
 have
\begin{align}
& \left \|\Box_k \mathscr{A} \partial^{3}_{x_i} f
\right\|_{L^\infty_{x_i} L^2_{(x_j)_{j\not=i}} L^2_t
(\mathbb{R}^{1+n})}  \lesssim  \|\Box_k f\|_{L^1_{x_i}
L^2_{(x_j)_{j\not=i}} L^2_t
(\mathbb{R}^{1+n}) },  \label{1-sm-mod} \\
& \left \|\Box_k \mathscr{A} \partial^{3}_{x_i} f
\right\|_{L^\infty_t L^2_x(\mathbb{R}^{1+n})} \lesssim  \langle
k_{\rm max}\rangle^{3/2} \| \Box_k f\|_{L^1_{x_i}
L^2_{(x_j)_{j\not=i}} L^2_t (\mathbb{R}^{1+n}) }. \label{1/2-sm-mod}
\end{align}
\end{prop}

{\bf Proof.} \eqref{1-sm-mod} holds by Proposition \ref{GSE1}
directly. In the case $|k_i| \ge 4$, \eqref{1/2-sm-mod} holds by
Proposition \ref{GSE3} and Lemma \ref{discret-deriv}. In the case
$|k_i| \le 3$, in view of Proposition \ref{GSE3},
\begin{align*}
& \left \|\Box_k \mathscr{A} \partial^{3}_{x_i} f
\right\|_{L^\infty_t L^2_x(\mathbb{R}^{1+n})} \lesssim  \left
\|D^{-3/2}_{x_i} \Box_k \mathscr{A} \partial^{3}_{x_i} f
\right\|_{L^\infty_t L^2_x(\mathbb{R}^{1+n})} \lesssim  \| \Box_k
f\|_{L^1_{x_i} L^2_{(x_j)_{j\not=i}} L^2_t (\mathbb{R}^{1+n}) },
\end{align*}
which implies the result, as desired. $\hfill \Box$

By the duality and Christ-Kiselev's Lemma in anisotropic Lebesgue
spaces \cite{WangNew}, we have the following
\begin{prop}\label{ma-sm-mo}
For $2\le  q \le \infty, q>8/n$ and  $k=(k_1,...,k_n)\in
 \mathbb{Z}^n, \abs{k_{i}} = k_{\rm max}$, $
i=1,...,n$,  we have
\begin{align}
 \left \|\Box_k \mathscr{A} \partial^{3}_{x_i} f \right\|_{L^{q}_{x_i}
L^\infty_{(x_j)_{j\not=i}} L^\infty_t (\mathbb{R}^{1+n}) } &
\lesssim \langle k_{\rm max}\rangle^{3/2+ 3/q}  \| \Box_k
f\|_{L^1_{x_i} L^2_{(x_j)_{j\not=i}} L^2_t (\mathbb{R}^{1+n}) }.
\label{ma-sm-mo-2}
\end{align}
\end{prop}

\begin{prop}\label{st-sm-mo}
 Let $2\le r < \infty$, $4/\gamma(r)=n(1/2-1/r)$ and $\gamma> \gamma(r)\vee 2$.  We have
\begin{align}
& \left \|\Box_k S(t) u_0 \right\|_{ L^\gamma_t
L^r_x(\mathbb{R}^{1+n}) } \lesssim   \|\Box_k u_0\|_{L^2
(\mathbb{R}^{n}) },  \label{st-sm-mo-1} \\
& \left \|\Box_k \mathscr{A}  f \right\|_{L^\infty_t L^2_x \, \cap
\, L^\gamma_t L^r_x (\mathbb{R}^{1+n}) } \lesssim    \|\Box_k f\|_{
L^{\gamma'}_t L^{r'}_x (\mathbb{R}^{1+n})
}, \label{st-sm-mo-2}\\
& \left \|\Box_k \mathscr{A} \partial^{3}_{x_i}  f
\right\|_{L^\gamma_t L^r_x (\mathbb{R}^{1+n}) } \lesssim \langle
k_{\rm max}\rangle^{3/2} \| \Box_k f\|_{L^1_{x_i}
L^2_{(x_j)_{j\not=i}} L^2_t
(\mathbb{R}^{1+n})}, \label{st-sm-mo-3}\\
& \left \|\Box_k \mathscr{A} \partial^{3}_{x_i}  f
\right\|_{L^\infty_{x_i} L^2_{(x_j)_{j\not=i}} L^2_t
(\mathbb{R}^{1+n})}   \lesssim \langle k_{\rm max}\rangle^{3/2}\|
\Box_k f\|_{ L^{\gamma'}_t L^{r'}_x (\mathbb{R}^{1+n})},
\label{st-sm-mo-4}
\end{align}
and for $2\le q<\infty$, $q>8/n$, $\alpha=0,3$,
\begin{align}
\left \|\Box_k \mathscr{A} \partial^\alpha_{x_i}  f
\right\|_{L^q_{x_i} L^\infty_{(x_j)_{j\not=i}} L^\infty_t
(\mathbb{R}^{1+n})} & \lesssim \langle k_{\rm max}\rangle^{\alpha+
3/q}
 \| \Box_k f\|_{ L^{\gamma'}_t L^{r'}_x
(\mathbb{R}^{1+n})}, \label{st-sm-mo-6}
\end{align}
\end{prop}
{\bf Proof.} \eqref{st-sm-mo-1} and \eqref{st-sm-mo-2} hold by
\ref{Strichartz-mod}. We now show \eqref{st-sm-mo-3}. We use the
same notations as in Proposition \ref{ma-sm-mo}. By Lemmas
\ref{Strichartz-mod}, \ref{discret-deriv} and Proposition
\ref{1order-sm-mod},
\begin{align}
 \mathcal {L}_k (\partial^{3}_{x_1} f, \psi) & \lesssim \langle k_{\rm max}\rangle^{3/2} \|\Box_k  f\|_{L^{1}_{x_1}
L^2_{x_2,...,x_n}L^2_t(\mathbb{R}^{1+n})} \left\| \tilde{\Box}_k
\psi \right\|_{L^{\gamma'}_t L^{r'}_{x} (\mathbb{R}^{1+n})}
\nonumber\\
& \lesssim \langle k_{\rm max}\rangle^{3/2} \|\Box_k
f\|_{L^{1}_{x_1} L^2_{x_2,...,x_n}L^2_t(\mathbb{R}^{1+n})} \left\|
\psi \right\|_{L^{\gamma'}_t L^{r'}_{x} (\mathbb{R}^{1+n})}.
 \label{st-sm-mo-7}
\end{align}
By duality,  \eqref{st-sm-mo-3}  holds by \eqref{st-sm-mo-7} and
Christ-Kiselev's Lemma. Exchanging the roles of $f$ and $\psi$, we
immediately have \eqref{st-sm-mo-4}. \eqref{st-sm-mo-6} holds by
Lemmas \ref{Strichartz-mod}, \ref{discret-deriv} and Proposition
\ref{MaxFunct-Mod-Dual}. $\hfill\Box$

we summarize the main conclusion of this section as following:
\begin{cor}\label{st-sm-m-c}
For $8/n \le p < \infty$, $2\le q < \infty$, $q>8/n$ and
 $k=(k_1,...,k_n)\in
 \mathbb{Z}^n$,
 $ \abs{k_{i}} = k_{\rm max}, i=1,...n$. We have
\begin{align}
& \left \| D^{3/2}_{x_i} \Box_k S(t) u_0 \right\|_{L^\infty_{x_i}
L^2_{(x_{j})j\neq i} L^2_t (\mathbb{R}^{1+n})}   \lesssim \| \Box_k
u_0\|_{ L^2 (\mathbb{R}^{n})}, \label{st-sm-m-c2}\\
& \left \| \Box_k S(t) u_0 \right\|_{L^q_{x_i}
L^\infty_{(x_{j})j\neq i} L^\infty_t (\mathbb{R}^{1+n})}   \lesssim
\langle k_{\rm max} \rangle^{3/q} \| \Box_k u_0\|_{ L^2
(\mathbb{R}^{n})},
\label{st-sm-m-c3} \\
 & \left \|\Box_k S(t) u_0 \right\|_{
L^{2+p}_{t,x}\, \cap \, L^\infty_t L^2_x (\mathbb{R}^{1+n}) }
\lesssim \|\Box_k u_0\|_{L^2 (\mathbb{R}^{n}) }, \label{st-sm-m-c1}
\end{align}
\begin{align}
& \left \|\Box_k \mathscr{A} \partial^{3}_{x_i} f
\right\|_{L^\infty_{x_i} L^2_{(x_{j})j\neq i} L^2_t
(\mathbb{R}^{1+n})} \lesssim  \| \Box_k f\|_{L^1_{x_i}
L^2_{(x_{j})j\neq i} L^2_t (\mathbb{R}^{1+n})},
\label{st-sm-mo-5}\\
& \left \|\Box_k \mathscr{A} \partial^{3}_{x_i}  f
\right\|_{L^q_{x_i} L^\infty_{(x_{j})j\neq i} L^\infty_t
(\mathbb{R}^{1+n})}  \lesssim \langle k_{\rm max}\rangle^{3/2+3/q}
\| \Box_k f\|_{L^1_{x_i} L^2_{(x_{j})j\neq i} L^2_t
(\mathbb{R}^{1+n})},
\label{st-sm-m-c6}\\
& \left \|\Box_k \mathscr{A} \partial^{3}_{x_i} f
\right\|_{L^\infty_t L^2_x \, \cap \, L^{2+p}_{t,x}
(\mathbb{R}^{1+n}) } \lesssim  \langle k_{\rm max}\rangle^{3/2}
\|\Box_k f\|_{ L^1_{x_i} L^2_{(x_{j})j\neq i} L^2_t
(\mathbb{R}^{1+n}) }. \label{st-sm-m-c4}
\end{align}
\begin{align}
&\left \|\Box_k \mathscr{A} \partial^{3}_{x_i} f
\right\|_{L^\infty_{x_i} L^2_{x_2,...,x_n} L^2_t (\mathbb{R}^{1+n})}
\lesssim \langle k_{\rm max}\rangle^{3/2}\| \Box_k f\|_{
L^{(2+p)/(1+p)}_{t,x} (\mathbb{R}^{1+n})},
\label{st-sm-mo-7}\\
& \left \|\Box_k \mathscr{A} \partial^{3}_{x_i}  f
\right\|_{L^q_{x_i} L^\infty_{x_2,...,x_n} L^\infty_t
(\mathbb{R}^{1+n})}  \lesssim
 \langle k_{\rm max} \rangle^{3+3/q} \| \Box_k
f\|_{L^{(2+p)/(1+p)}_{t,x} (\mathbb{R}^{1+n})},
\label{st-sm-m-c8}\\
& \left \|\Box_k \mathscr{A}  f \right\|_{L^\infty_t L^2_x \, \cap
\, L^{2+p}_{t,x}  (\mathbb{R}^{1+n}) } \lesssim    \|\Box_k f\|_{
L^{(2+p)/(1+p)}_{t,x} (\mathbb{R}^{1+n}) }. \label{st-sm-m-c9}
\end{align}
\end{cor}

\section{Linear estimates with derivative interaction} \label{LE2}

Recall that in Prop \ref{GSE1} we assume that $\abs{k_{i}} = k_{\rm
max}$ for any $k\in\Z^n$.  In view of \eqref{st-sm-mo-5} in
Corollary \ref{st-sm-m-c}, the operator $\mathscr{A}$ in the space
$L^\infty_{x_1} L^2_{x_2,...,x_n} L^2_t (\mathbb{R}^{1+n})$ has
succeed in absorbing the partial derivative $\partial^{3}_{x_1}$.
However,  it seem that $\mathscr{A}$ can not deal with the partial
derivative $\partial^{3}_{x_2}$ in the space
 $L^\infty_{x_1} L^2_{x_2,...,x_n} L^2_t (\mathbb{R}^{1+n})$. So, we need a new way to handle the interaction
between $L^\infty_{x_1} L^2_{x_2,...,x_n} L^2_t (\mathbb{R}^{1+n})$
and $\partial^{3}_{x_2}$.

\begin{prop}\label{sm-eff-interact1}
For $i=2,...,n$, $2\le q \le \infty$, $q>8/n$. Let $4\le r <\infty$,
$2/\gamma(r)= n(1/2-1/r)$, $\gamma> 2\vee \gamma(r)
4$.\\
For $\abs{k_{1}}= k_{\rm max}$,  we have
 \begin{align}
& \left \|\Box_k \partial^{3}_{x_i}  \mathscr{A}  f
\right\|_{L^\infty_{x_1} L^2_{x_2,...,x_n} L^2_t (\mathbb{R}^{1+n})}
\lesssim  \| \partial^{3}_{x_i} \partial^{-3}_{x_1} \Box_k
f\|_{L^1_{x_1} L^2_{x_2,...,x_n}L^2_t (\mathbb{R}^{1+n})},
\label{sm-int-1}\\
& \left \|\Box_k \partial^{3}_{x_i}  \mathscr{A}  f
\right\|_{L^\infty_{x_1} L^2_{x_2,...,x_n} L^2_t (\mathbb{R}^{1+n})}
\lesssim  \| \partial^{3}_{x_i} D^{-3/2}_{x_1} \Box_k
f\|_{L^{\gamma'} L^{r'}_x (\mathbb{R}^{1+n})}, \label{sm-int-1a}
\end{align}
For $\abs{k_{i}} = k_{\rm max}$, we have
\begin{align}
& \left \|\Box_k \partial^{3}_{x_i}  \mathscr{A}  f
\right\|_{L^q_{x_1} L^\infty_{x_2,...,x_n} L^\infty_t
(\mathbb{R}^{1+n})} \lesssim \langle k_{\rm
max}\rangle^{\frac{3}{2}+\frac{3}{q}}   \|
 \Box_k  f\|_{L^1_{x_i}
L^2_{(x_j)_{j\not=i}}L^2_t (\mathbb{R}^{1+n})},
\label{sm-int-2}\\
& \left \|\Box_k \partial^{3}_{x_i}  \mathscr{A}  f
\right\|_{L^q_{x_1} L^\infty_{x_2,...,x_n} L^\infty_t
(\mathbb{R}^{1+n})} \lesssim \langle k_{\rm
max}\rangle^{3+\frac{3}{q}} \|
 \Box_k  f\|_{L^{\gamma'}_{t}
L^{r'}_{x}  (\mathbb{R}^{1+n})}. \label{sm-int-2a}
\end{align}
\end{prop}

{\bf Proof.} From Proposition \ref{GSE1}, we can get
\eqref{sm-int-1} directly. As before, we only give the proof when
$i=2$. For \eqref{sm-int-1a}, because of
\begin{align}
 \mathcal {L} (\partial^{3}_{x_2} \Box_k  f, \psi) &:= \left|\int_\mathbb{R}\left(
\int_{\mathbb{R}}S(t-\tau) \partial^{3}_{x_2} \Box_k
f(\tau) d\tau, \  \psi(t) \right ) dt \right| \nonumber\\
& \le \left\| \int_{\mathbb{R}}S(-\tau) \partial^{3}_{x_2}
D^{-3/2}_{x_1} \Box_k  f(\tau) d\tau\right\|_{L^2(\mathbb{R}^n)}
\left\|\tilde{\Box}_k D^{3/2}_{x_1} \int_\mathbb{R} S(-t) \psi(t) dt
\right\|_{L^2(\mathbb{R}^n)}.
 \label{sm-int-3}
\end{align}
By the Strichartz inequality and Proposition \ref{GSE3},
\begin{align}
 \mathcal {L} (\partial^{3}_{x_2} f, \psi)
& \lesssim \| \partial^{3}_{x_2} D^{-3/2}_{x_1} \Box_k
f\|_{L^{\gamma'}_t L^{r'}_x (\mathbb{R}^{1+n})}  \|\psi
\|_{L^1_{x_1} L^2_{x_2,...,x_n}L^2_t (\mathbb{R}^{1+n})} .
 \label{sm-int-3a}
\end{align}
By duality, \eqref{sm-int-3a} implies \eqref{sm-int-1a}. The proof
of \eqref{sm-int-2} is similar. From Propositions \ref{GSE3},
\ref{MaxFunct-Mod-Dual} and Lemma \ref{discret-deriv},
\begin{align}
 \mathcal {L} (\partial^{3}_{x_2} \Box_k f, \psi)
& \le \left\| \int_{\mathbb{R}}S(-\tau) D^3_{x_2} \Box_k f(\tau)
d\tau\right\|_{L^2(\mathbb{R}^n)} \left\|\tilde{\Box}_k
\int_\mathbb{R} S(-t) \psi(t) dt \right\|_{L^2(\mathbb{R}^n)}
\nonumber\\
& \lesssim \langle k_{\rm max}\rangle^{3/2} \|  \Box_k
f\|_{L^1_{x_2} L^2_{(x_j)_{j\not=2}}L^2_t (\mathbb{R}^{1+n})}
\langle k_{\rm max} \rangle^{3/q} \|\tilde{\Box}_k \psi
\|_{L^{q'}_{x_1}
L^1_{x_2,...,x_n}L^1_t (\mathbb{R}^{1+n})} \nonumber\\
& \lesssim \langle k_{\rm max}\rangle^{3/2}\langle k_{\rm max}
\rangle^{3/q} \| \Box_k f\|_{L^1_{x_2} L^2_{(x_j)_{j\not=2}}L^2_t
(\mathbb{R}^{1+n})} \|\psi \|_{L^{q'}_{x_1} L^1_{x_2,...,x_n}L^1_t
(\mathbb{R}^{1+n})}.
 \label{sm-int-4}
\end{align}
\eqref{sm-int-2} is the duality of \eqref{sm-int-4}. For
\eqref{sm-int-2a}, noticing that
\begin{align}
 \mathcal {L} (\partial^{3}_{x_2} \Box_k f, \psi)
& \le \left\| \int_{\mathbb{R}}S(-\tau) \partial^{3}_{x_2} \Box_k
f(\tau) d\tau\right\|_{L^2(\mathbb{R}^n)} \left\|\tilde{\Box}_k
\int_\mathbb{R} S(-t) \psi(t) dt \right\|_{L^2(\mathbb{R}^n)}
\nonumber\\
& \lesssim \langle k_{\rm max}\rangle^{3}  \langle k_{\rm
max}\rangle^{3/q} \|\psi\|_{L^{q'}_{x_1} L^1_{x_2,...,x_n}L^1_t
(\mathbb{R}^{1+n})} \|
 \Box_k
f \|_{L^{\gamma'}_{t} L^{r'}_{x} (\mathbb{R}^{1+n})} ,
 \label{sm-int-4a}
\end{align}
which implies \eqref{sm-int-2a}, as desired.  $\hfill\Box$

\begin{lem}\label{sm-ef-int2}
Let  $\psi : [0, \infty) \to [0,1]$ be a smooth bump function
satisfying $\psi(x)=1$ as $|x| \le 1$ and $\psi (x)=0$ if $|x| \ge
2$. Denote $\psi_1(\xi) = \psi(\xi_2/2\xi_1)$,  $\psi_2(\xi) =
1-\psi(\xi_2/2\xi_1)$, $\xi \in \mathbb{R}^n$. Then  we have for
$\sigma \ge 0$,
 \begin{align}
& \sum_{k\in \mathbb{Z}^n, \, |k_1|>4} \langle k_1\rangle^\sigma
\left \|\mathscr{F}^{-1}_{\xi_1, \xi_2} \psi_1 \mathscr{F}_{x_1,
x_2} \Box_k
\partial^{3}_{x_2} \mathscr{A} f \right\|_{L^\infty_{x_1}
L^2_{x_2,...,x_n} L^2_t (\mathbb{R}^{1+n})}  \nonumber\\
& \quad \quad\quad  \lesssim \sum_{k\in \mathbb{Z}^n, \, |k_1|>4}
\langle k_1\rangle^\sigma \left \|\Box_k
 f \right\|_{L^1_{x_1}
L^2_{x_2,...,x_n} L^2_t (\mathbb{R}^{1+n})},  \label{sm-int-5a}
\end{align}
and for $\sigma \ge 3$,
 \begin{align}
& \sum_{k\in \mathbb{Z}^n, \, |k_1|>4} \langle k_1\rangle^\sigma
\left \|\mathscr{F}^{-1}_{\xi_1, \xi_2} \psi_2 \mathscr{F}_{x_1,
x_2} \Box_k
\partial^{3}_{x_2} \mathscr{A} f \right\|_{L^\infty_{x_1}
L^2_{x_2,...,x_n} L^2_t (\mathbb{R}^{1+n})}  \nonumber\\
 &  \quad \quad\quad \lesssim \sum_{k\in \mathbb{Z}^n, \,
|k_2|>4} \langle k_1\rangle^{\sigma-3}\langle k_2\rangle^{3} \left
\|\Box_k
 f \right\|_{L^1_{x_1}
L^2_{x_2,...,x_n} L^2_t (\mathbb{R}^{1+n})}. \label{sm-int-5b}
\end{align}
\end{lem}

{\bf Proof.}  For the terseness of proof, we let
$$
I= \left \|\mathscr{F}^{-1}_{\xi_1, \xi_2} \psi_1 \mathscr{F}_{x_1,
x_2} \Box_k
\partial^{3}_{x_2} \mathscr{A} f \right\|_{L^\infty_{x_1}
L^2_{x_2,...,x_n} L^2_t (\mathbb{R}^{1+n})},
 $$
$$
II= \left \|\mathscr{F}^{-1}_{\xi_1, \xi_2} \psi_2 \mathscr{F}_{x_1,
x_2} \Box_k
\partial^{3}_{x_2} \mathscr{A} f \right\|_{L^\infty_{x_1}
L^2_{x_2,...,x_n} L^2_t (\mathbb{R}^{1+n})}.
$$

 Firstly, we give the estimate of $I$. Let $\eta_k $ be as in Lemma \ref{eq-mod-funct}. For $k \in
\mathbb{Z}^n$, $|k_1|
> 4$, applying the almost orthogonality of $\Box_k$,  we have
\begin{align}
I & \lesssim  \sum_{|\ell_1|, |\ell_2| \le 1}  \left
\|\mathscr{F}^{-1}_{\xi_1, \xi_2} \psi
\left(\frac{\xi_2}{2\xi_1}\right) (\frac{\xi_2}{\xi_1})^{3}
\prod_{i=1,2}\eta_{k_i+\ell_i} (\xi_i) \mathscr{F}_{x_1, x_2} \Box_k
\partial^{3}_{x_1} \mathscr{A} f \right\|_{L^\infty_{x_1}
L^2_{x_2,...,x_n} L^2_t (\mathbb{R}^{1+n})}. \label{sm-int-7}
\end{align}
Denote
\begin{align}
(f \circledast_{12} g)(x) = \int_{\mathbb{R}^2}  f(t, x_1-y_1,
x_2-y_2, x_3,..., x_n) g(t, y_1,y_2) dy_1 dy_2. \label{sm-int-8}
\end{align}
We have for any Banach function space $X$ defined on
$\mathbb{R}^{1+n}$,
\begin{align}
\|f \circledast_{12} g\|_X \le \|g\|_{L^1_{y_1,y_2}(\mathbb{R}^2)}
\sup_{y_1, y_2} \| f(\cdot, \cdot-y_1, \cdot-y_2, \cdot,...,
\cdot)\|_X. \label{sm-int-9}
\end{align}
Hence, combining \eqref{sm-int-7} with \eqref{sm-int-9},
\begin{align}
I & \lesssim  \sum_{|\ell_1|, |\ell_2| \le 1}  \left
\|\mathscr{F}^{-1}_{\xi_1, \xi_2} \psi
\left(\frac{\xi_2}{2\xi_1}\right) (\frac{\xi_2}{\xi_1})^{3}
\prod_{i=1,2}\eta_{k_i+\ell_i} (\xi_i) \right\|_{L^1(\mathbb{R}^2)}
\left\|\Box_k
\partial^{3}_{x_1} \mathscr{A} f \right\|_{L^\infty_{x_1}
L^2_{x_2,...,x_n} L^2_t (\mathbb{R}^{1+n})}. \label{sm-int-10}
\end{align}
Using Bernstein's multiplier estimate, for $|k_1|>4$, we have
\begin{align}
& \!\! \!\! \left \|\mathscr{F}^{-1}_{\xi_1, \xi_2} \psi
\left(\frac{\xi_2}{2\xi_1}\right) (\frac{\xi_2}{\xi_1})^{3}
\prod_{i=1,2}\eta_{k_i+\ell_i} (\xi_i) \right\|_{L^1(\mathbb{R}^2)}
\nonumber\\
& \lesssim  \sum_{|\alpha| \le 2} \left \|D^\alpha \left[\psi
\left(\frac{\xi_2}{2\xi_1}\right) (\frac{\xi_2}{\xi_1})^{3}
\prod_{i=1,2}\eta_{k_i+\ell_i} (\xi_i) \right]
\right\|_{L^2(\mathbb{R}^2)} \lesssim 1.
 \label{sm-int-11}
\end{align}
By Proposition \ref{1order-sm-mod}, \eqref{sm-int-10} and
\eqref{sm-int-11}, we have
\begin{align}
I & \lesssim \left\|\Box_k f \right\|_{L^1_{x_1} L^2_{x_2,...,x_n}
L^2_t (\mathbb{R}^{1+n})}, \ \ |k_1| >4. \label{sm-int-12}
\end{align}
Next, we consider the estimate of $II$. Using Proposition
\ref{sm-eff-interact1},
\begin{align}
II & \lesssim \left \| \mathscr{F}^{-1}_{\xi_1, \xi_2}(\xi_2/\xi_1)
\psi_2 \mathscr{F}_{x_1,x_2}  \Box_k
 f \right\|_{L^1_{x_1}
L^2_{x_2,...,x_n} L^2_t (\mathbb{R}^{1+n})} \nonumber\\
& \lesssim \sum_{|\ell_1|, |\ell_2| \le 1}  \left
\|\mathscr{F}^{-1}_{\xi_1, \xi_2} \left(1-\psi \left(\frac{\xi_2}{2
\xi_1}\right)\right) (\frac{\xi_2}{\xi_1})^{3}
\prod_{i=1,2}\eta_{k_i+\ell_i} (\xi_i) \right\|_{L^1(\mathbb{R}^2)}
\nonumber\\
& \quad \quad \quad \times  \left\|\Box_k
 f \right\|_{L^1_{x_1}
L^2_{x_2,...,x_n} L^2_t (\mathbb{R}^{1+n})}. \label{sm-int-13}
\end{align}
Notice that ${\rm supp }\psi_2 \subset \{\xi: \, |\xi_2| \ge
2|\xi_1|\}$. If $|k_1| > 4$, we have $|k_2|>6$ and $|k_2| \ge |k_1|$
in the summation of the left-hand side of \eqref{sm-int-5b}. So,
$\sum_{k\in \mathbb{Z}^n, \ |k_1|>4}\langle k_1\rangle^{\sigma} II
\le \sum_{k\in \mathbb{Z}^n, \ |k_1|>4}\langle k_2\rangle^{\sigma-3}
\langle k_1\rangle^{3}  II$.
\begin{align}
& \left \|\mathscr{F}^{-1}_{\xi_1, \xi_2} \left(1-\psi
\left(\frac{\xi_2}{2 \xi_1}\right)\right) (\frac{\xi_2}{\xi_1})^{3}
\prod_{i=1,2}\eta_{k_i+\ell_i} (\xi_i) \right\|_{L^1(\mathbb{R}^2)} \nonumber\\
& \lesssim  \sum_{|\alpha|\le 2} \left \|D^\alpha
\left[\mathscr{F}^{-1}_{\xi_1, \xi_2} \left(1-\psi
\left(\frac{\xi_2}{2 \xi_1}\right)\right) (\frac{\xi_2}{\xi_1})^{3}
\prod_{i=1,2}\eta_{k_i+\ell_i} (\xi_i)\right]
\right\|_{L^2(\mathbb{R}^2)} \nonumber\\
&  \lesssim \langle k_2 \rangle^{3}  \langle k_1\rangle^{-3}.
 \label{sm-int-14}
\end{align}
Collecting \eqref{sm-int-12} and \eqref{sm-int-14}, we get the
result, as desired. $\hfill \Box$

\begin{lem}\label{sm-ef-int4}
Let $k=(k_1,...,k_n)$, $2\le q \le \infty$, $q>8/n$. Then we have
for $\sigma \ge 0$ and $i, \alpha=1,...,n$,
 \begin{align}
& \sum_{k\in \mathbb{Z}^n, \, |k_\alpha|=k_{\rm max}>4} \langle k
\rangle^\sigma \left \|\Box_k
\partial^{3}_{x_i} \mathscr{A} f \right\|_{L^q_{x_1}
L^\infty_{x_2,...,x_n} L^\infty_t (\mathbb{R}^{1+n})}  \nonumber\\
 &  \quad \quad\quad \lesssim \sum_{k\in \mathbb{Z}^n, \,
|k_\alpha|=k_{\rm max}>4} \langle k_\alpha \rangle^{\sigma+ 3/2+3/q}
\left \|\Box_k
 f \right\|_{L^1_{x_\alpha}
L^2_{(x_j)_{j\not=\alpha}} L^2_t (\mathbb{R}^{1+n})}.
\label{sm-int-4.4a}
\end{align}
\end{lem}

{\bf Proof.} First, we consider the case $\alpha=1$. In view of
\eqref{st-sm-m-c6} and $|k_1|= k_{\rm max}>4$,
\begin{align}
  \left \|\Box_k
\partial^{3}_{x_i} \mathscr{A} f \right\|_{L^{q}_{x_1}
L^\infty_{x_2,...,x_n} L^\infty_t (\mathbb{R}^{1+n})}  & \ \lesssim
\sum_{|\ell_1|, |\ell_i| \le 1} \left\| \mathscr{F}^{-1}_{\xi_1,
\xi_i}  \left((\frac{\xi_i}{\xi_1})^{3} \eta_{k_i+\ell_i} (\xi_i)
\eta_{k_1+\ell_1} (\xi_1) \right)\right\|_{L^1(\mathbb{R}^2)}
\nonumber\\
 & \ \ \ \   \times \| \Box_k
\partial^{3}_{x_1} \mathscr{A} f \|_{L^{q}_{x_1}
L^\infty_{x_2,...,x_n} L^\infty_t (\mathbb{R}^{1+n})} \nonumber\\
 & \ \lesssim \langle k_i\rangle^{3} \langle k_1\rangle^{-3} \langle k_1\rangle^{3/2+3/q} \| \Box_k
 f \|_{L^{1}_{x_1} L^2_{x_2,...,x_n} L^2_t
(\mathbb{R}^{1+n})} \nonumber\\
 & \ \lesssim  \langle k_1\rangle^{3/2+3/q} \| \Box_k
 f \|_{L^{1}_{x_1} L^2_{x_2,...,x_n} L^2_t
(\mathbb{R}^{1+n})} . \label{sm-int-4.4b}
\end{align}
\eqref{sm-int-4.4b} implies the result, as desired. Next, we
consider the case $\alpha=2$, the general case $\alpha=3,...,n$ is
similar. Notice that $|k_2| = k_{\rm max} >4$. By \eqref{sm-int-2},
\begin{align}
  \left \|\Box_k
\partial^{3}_{x_i} \mathscr{A} f \right\|_{L^{q}_{x_1}
L^\infty_{x_2,...,x_n} L^\infty_t (\mathbb{R}^{1+n})}  & \ \lesssim
\sum_{|\ell_2|, |\ell_i| \le 1} \left\| \mathscr{F}^{-1}_{\xi_i,
\xi_2}  \left((\frac{\xi_i}{\xi_2})^{3} \eta_{k_i+\ell_i} (\xi_i)
\eta_{k_2+\ell_2} (\xi_2) \right)\right\|_{L^1(\mathbb{R}^2)}
\nonumber\\
 & \ \ \ \   \times \| \Box_k
\partial^{3}_{x_2} \mathscr{A} f \|_{L^{q}_{x_1}
L^\infty_{x_2,...,x_n} L^\infty_t (\mathbb{R}^{1+n})} \nonumber\\
 & \ \lesssim \langle k_i\rangle^{3} \langle k_2\rangle^{-3} \langle k_2\rangle^{3/2+3/q} \| \Box_k
 f \|_{L^{1}_{x_2} L^2_{x_1,x_3,...,x_n} L^2_t
(\mathbb{R}^{1+n})} \nonumber\\
 & \ \lesssim  \langle k_2\rangle^{3/2+3/q} \| \Box_k
 f \|_{L^{1}_{x_2} L^2_{x_1,x_3,...,x_n} L^2_t
(\mathbb{R}^{1+n})} . \label{sm-int-4.4c}
\end{align}
$\hfill \Box$

\begin{rem} \rm
From the proof of Lemma \ref{sm-ef-int4}, we easily see that for
$\beta=1,...,n$
\begin{align}
& \sum_{k\in \mathbb{Z}^n, \, |k_\alpha|=k_{\rm max}>4} \langle k
\rangle^\sigma \left \|\Box_k
\partial^{3}_{x_i} \mathscr{A} f \right\|_{L^q_{x_\beta}
L^\infty_{(x_j)_{j\not=\beta}} L^\infty_t (\mathbb{R}^{1+n})}  \nonumber\\
 &  \quad \quad\quad \lesssim \sum_{k\in \mathbb{Z}^n, \,
|k_\alpha|=k_{\rm max}>4} \langle k_\alpha \rangle^{\sigma+ 3/2+3/q}
\left \|\Box_k
 f \right\|_{L^1_{x_\alpha}
L^2_{(x_j)_{j\not=\alpha}} L^2_t (\mathbb{R}^{1+n})}.
\label{sm-int-4.4d}
\end{align}
\end{rem}

\section{Proof of Theorem \ref{DNLS1-mod}} \label{pf-thm1}
Before we prove our main Theorem \ref{DNLS-mod}, we consider it's
special version-Theorem \ref{DNLS1-mod}. As assumption, the
nonlinear term takes the form
\begin{align*}
F((\partial_x^\alpha u)_{\abs{\alpha}\leq 3},
(\partial_x^\alpha\bar{u})_{\abs{\alpha}\leq 3}) =
\sum_{i=1}^{n}\lambda_i\partial^{3}_{x_i} (u^{\kappa_i+1}).
\end{align*}

Denote
\begin{align*}
& \rho_1(u) = \sum^n_{i=1}\sum_{k\in \mathbb{Z}^n, \ |k_i|= k_{\rm
max} >4} \langle k_i\rangle^{3}  \|\Box_k u\|_{L^{\infty}_{x_i}
L^2_{(x_j)_{j\not=i}}L^2_t(\mathbb{R}^{1+n})},\\
& \rho_2(u) = \sum^n_{i=1}\sum_{k\in \mathbb{Z}^n} \langle
k\rangle^{3/2-3/\kappa} \|\Box_k u\|_{L^{\kappa}_{x_i}
L^\infty_{(x_j)_{j\not=i}}L^\infty_t(\mathbb{R}^{1+n})},\\
& \rho_3(u) = \sum_{k\in \mathbb{Z}^n} \langle k\rangle^{3/2}
\|\Box_k u\|_{L^\infty_tL^2_x  \bigcap L^{2+\kappa}_{x, t}
(\mathbb{R}^{1+n})}.
\end{align*}

Put
\begin{align*}
X:= \left\{u\in \mathscr{S}'(\mathbb{R}^{1+n}): \ \|u\|_X:=
\sum^3_{i=1}\rho_i(u) \le \delta_0 \right\}.
\end{align*}

Considering the following integral mapping:
\begin{align*}
\mathscr{T}: u(t) \to S(t) u_0 - {\rm i} \mathscr{A}
\left(\sum^n_{i=1} \lambda_i \partial^{3}_{x_i} u^{\kappa_i+1}
\right).
\end{align*}

Firstly, we estimate  $\ \|S(t)u_0\|_{X}$. \\
For simplicity, we denote
\begin{align*}
& \|u\|_{Y_i} = \sum_{k\in \mathbb{Z}^n, \ |k_i|= k_{\rm max} >4}
\langle k_i\rangle^{3}  \|\Box_k u\|_{L^{\infty}_{x_i}
L^2_{(x_j)_{j\not=i}}L^2_t(\mathbb{R}^{1+n})}.
\end{align*}
For the estimate of $\rho_1(u)$, it suffices to control
$\|\cdot\|_{Y_1}$. By \eqref{smo-eff.21} and Plancherel's identity,
we have
\begin{align*}
\|S(t) u_0\|_{Y_1} & \lesssim  \sum_{k\in \mathbb{Z}^n, \ |k_1|=
k_{\rm max} >4} \langle k_1\rangle^{3}  \|\Box_k D^{-3/2}_{x_1}
u_0\|_{L^2(\mathbb{R}^{n})} \\
&  \lesssim \sum_{k\in \mathbb{Z}^n} \langle k_1\rangle^{3/2}
\|\Box_k u_0\|_{L^2(\mathbb{R}^{n})}.
\end{align*}
For  $\rho_2(u)$, using Prop \ref{MaxFunct-Mod}, we obtain
\begin{align*}
 \|\Box_k S(t)
u_0\|_{L^{q}_{x_i} L^\infty_{(x_j)_{j\not=i}} L^\infty_t
(\mathbb{R}^{1+n}) } & \lesssim  \langle k_{\rm max} \rangle^{3/q}
\| \Box_k u_0\|_{L^2(\mathbb{R}^n)}
\end{align*}
Therefore, we have
\begin{align*}
\rho_2 (S(t) u_0) &  \lesssim \sum_{k\in \mathbb{Z}^n} \langle k
\rangle^{3/2} \|\Box_k u_0\|_{L^2(\mathbb{R}^{n})}
\end{align*}
For  $\rho_3(u)$, using Lemma \ref{Strichartz-mod}, we have
\begin{align*}
\rho_3 (S(t) u_0) &  \lesssim \sum_{k\in \mathbb{Z}^n} \langle k
\rangle^{3/2} \|\Box_k u_0\|_{L^2(\mathbb{R}^{n})}
\end{align*}
Secondly, we estimate $\ \|\mathscr{A} \left(\sum^n_{i=1} \lambda_i
\partial^{3}_{x_i} u^{\kappa_i+1}\right)\|_{X}$.\\
For the simplicity of proof, we denote
\begin{align*}
& \mathbb{S}^{(i)}_{\ell, 1}:= \{k^{(1)},..., k^{(\kappa_\ell+1)}
\in \mathbb{Z}^n : \
|k^{(1)}_i|\vee...\vee |k^{(\kappa_\ell+1)}_i|>4\}, \\
& \mathbb{S}^{(i)}_{\ell, 2}:= \{k^{(1)},..., k^{(\kappa_\ell+1)}
\in \mathbb{Z}^n: \ |k^{(1)}_i|\vee...\vee
|k^{(\kappa_\ell+1)}_i|\le 4\}.
\end{align*}
Using the frequency-uniform decomposition, we have
\begin{align}
u^{\kappa_\ell+1}  =  & \sum_{k^{(1)},..., k^{(\kappa_\ell+1)}\in
\mathbb{Z}^n} \Box_{k^{(1)}} u ... \Box_{k^{(\kappa_\ell+1)}} u  \nonumber\\
 = & \sum_{\mathbb{S}^{(i)}_{\ell,1}} \Box_{k^{(1)}} u ... \Box_{k^{(\kappa_\ell+1)}} u  +
\sum_{\mathbb{S}^{(i)}_{\ell,2}} \Box_{k^{(1)}} u ...
\Box_{k^{(\kappa_\ell+1)}} u . \label{decomp}
\end{align}

Using \eqref{st-sm-mo-5} and \eqref{st-sm-mo-7}, we obtain that
\begin{align}
& \|\mathscr{A} \partial^{3}_{x_1} u^{\kappa_1+1}\|_{Y_1}  \lesssim
\sum_{k \in \mathbb{Z}^n, \ |k_1|= k_{\rm max} >4} \langle
k_1\rangle^{3} \sum_{\mathbb{S}^{(1)}_{1,1}} \|\Box_k
\left(\Box_{k^{(1)}} u ... \Box_{k^{(\kappa_1+1)}} u \right)
\|_{L^{1}_{x_1}
L^2_{x_2,...,x_n}L^2_t(\mathbb{R}^{1+n})}  \nonumber\\
& \ \  +  \sum_{k \in \mathbb{Z}^n, \ |k_1|= k_{\rm max} >4} \langle
k_1\rangle^{3}\langle k_{\rm max}\rangle^{3/2}
\sum_{\mathbb{S}^{(1)}_{1,2}} \| \Box_k \left(\Box_{k^{(1)}} u ...
\Box_{k^{(\kappa_1+1)}} u
\right)\|_{L^{(2+\kappa)/(1+\kappa)}_{t,x}(\mathbb{R}^{1+n}) } \nonumber\\
& := I +II . \label{nonhomo-1}
\end{align}
In view of the support property of $\widehat{\Box_k u}$, we see that
\begin{align}
 \Box_k \left(\Box_{k^{(1)}} u ... \Box_{k^{(\kappa_1+1)}} u
 \right)=0, \ \ if  \ \ |k- k^{(1)}-...-k^{(\kappa_1+1)}| \ge C.
 \label{orth}
\end{align}
Hence, by Lemma \ref{discret-deriv},
\begin{align}
I & \lesssim \sum_{k \in \mathbb{Z}^n, \ |k_1|= k_{\rm max} >4}
\langle k_1\rangle^{3}   \sum_{\mathbb{S}^{(1)}_{1,1}}
\|\Box_{k^{(1)}} u ... \Box_{k^{(\kappa_1+1)}} u  \|_{L^{1}_{x_1}
L^2_{x_2,...,x_n}L^2_t(\mathbb{R}^{1+n})} \chi_{|k-
k^{(1)}-...-k^{(\kappa_1+1)}| \le C}. \label{non-est-1}
\end{align}
By H\"older's inequality and $\|\Box_k u \|_{L^{\infty}_{x}}
\lesssim \|\Box_k u \|_{L^{2}_{x}} $ uniformly holds for all $k\in
\mathbb{Z}^n$, we have
\begin{align*}
& \|\Box_{k^{(1)}} u ... \Box_{k^{(\kappa_1+1)}} u \|_{L^{1}_{x_1}
L^2_{x_2,...,x_n}L^2_t(\mathbb{R}^{1+n})} \\
& \le \|\Box_{k^{(1)}} u\|_{L^{\infty}_{x_1}
L^2_{x_2,...,x_n}L^2_t(\mathbb{R}^{1+n})} \prod^{\kappa_1+1}_{i=2}
\|\Box_{k^{(i)}} u \|_{L^{\kappa}_{x_1}
L^\infty_{x_2,...,x_n}L^\infty_t \, \bigcap \, L^\infty_t
L^2_x(\mathbb{R}^{1+n})}.
\end{align*}
Since $|k- k^{(1)}-...-k^{(\kappa_1+1)}| \le C$ implies that $|k_1-
k^{(1)}_1-...-k^{(\kappa_1+1)}_1| \le C$, we see that $|k_1| \sim \
\max_{i=1,...,\kappa_1+1} |k^{(i)}_1|$. Without loss of generality ,
 we may assume that $|k^{(1)}_1|= \max_{i=1,...,\kappa_1+1}
|k^{(i)}_1|$ in the summation $\sum_{\mathbb{S}^{(1)}_{1,1}}$ in
\eqref{non-est-1} above. Therefore, we have
\begin{align}
I & \lesssim    \sum_{k^{(1)} \in \mathbb{Z}^n, \ |k^{(1)}_1| \sim
k_{\rm max} >4} \langle k^{(1)}_1\rangle^{3}  \|\Box_{k^{(1)}}
u\|_{L^{\infty}_{x_1} L^2_{x_2,...,x_n}L^2_t(\mathbb{R}^{1+n})} \nonumber\\
& \ \ \times  \sum_{k^{(2)},...,k^{(\kappa_1+1)}\in \mathbb{Z}^n}
\prod^{\kappa_1+1}_{i=2} \|\Box_{k^{(i)}} u \|_{L^{\kappa}_{x_1}
L^\infty_{x_2,...,x_n}L^\infty_t \, \bigcap \, L^\infty_t
L^2_x (\mathbb{R}^{1+n})} \nonumber\\
& \lesssim \rho_1(u) (\rho_2(u)+\rho_3(u))^{\kappa_1} .
\label{est-I}
\end{align}
In view of \eqref{orth} we easily see that $|k_1| \le C$ in $II$ of
\eqref{nonhomo-1}. Hence,
\begin{align}
II & \lesssim \sum_{k \in \mathbb{Z}^n, \ |k_1|= k_{\rm max} >4}
\sum_{\mathbb{S}^{(1)}_{1,2}} \|\Box_{k^{(1)}} u ...
\Box_{k^{(\kappa_1+1)}} u  \|_{L^{(2+\kappa)/(1+\kappa)}_{x, t}
(\mathbb{R}^{1+n})}
\chi_{|k- k^{(1)}-...-k^{(\kappa_1+1)}| \le C} \nonumber\\
& \lesssim  \sum_{\mathbb{S}^{(1)}_{1,2}} \|\Box_{k^{(1)}} u ...
\Box_{k^{(\kappa_1+1)}} u \|_{L^{(2+\kappa)/(1+\kappa)}_{x, t}
(\mathbb{R}^{1+n})} \nonumber\\
& \lesssim  \sum_{\mathbb{S}^{(1)}_{1,2}}
\prod^{\kappa_1+1}_{i=1}\|\Box_{k^{(i)}} u \|_{L^{2+\kappa}_{x, t}
\, \bigcap \, L^\infty_t L^2_x (\mathbb{R}^{1+n})} \lesssim
\rho_3(u)^{1+\kappa_1}. \label{non-est-2}
\end{align}
Hence, we have
\begin{align}
\|\mathscr{A} \partial^{3}_{x_1} u^{\kappa_1+1}\|_{Y_1} & \lesssim
\rho_1(u) (\rho_2(u)+\rho_3(u))^{\kappa_1} + \rho_3(u)^{1+\kappa_1}.
\label{est-lam1}
\end{align}
Now, we turn to estimate $\|\mathscr{A} \partial^{3}_{x_2}
u^{\kappa_2+1}\|_{Y_1}$. Let $\psi_i$ be as in Lemma
\ref{sm-ef-int2}.  For convenience, we write
\begin{align}
P_i= \mathscr{F}^{-1}_{\xi_1, \xi_2} \psi_i \mathscr{F}_{x_1, x_2},
\ \ i=1,2. \label{Proj-i}
\end{align}
We have
 \begin{align}
\|\mathscr{A} \partial^{3}_{x_2} u^{\kappa_2+1}\|_{Y_1} \lesssim &
\left \|P_1
\partial^{3}_{x_2} \mathscr{A} u^{\kappa_2+1} \right\|_{Y_1} + \left \|P_2
\partial^{3}_{x_2} \mathscr{A} u^{\kappa_2+1} \right\|_{Y_1}:=III+IV. \label{est-int-1}
\end{align}
Using the decomposition \eqref{decomp},
\begin{align}
III & \le \Big \|P_1
\partial^{3}_{x_2} \mathscr{A} \sum_{\mathbb{S}^{(1)}_{2,1}} (\Box_{k^{(1)}} u ...
\Box_{k^{(\kappa_2+1)}} u )  \Big\|_{Y_1} \nonumber\\
& \quad  +  \Big \|P_1
\partial^{3}_{x_2} \mathscr{A} \sum_{\mathbb{S}^{(1)}_{2,2}} (\Box_{k^{(1)}} u ...
\Box_{k^{(\kappa_2+1)}} u )  \Big\|_{Y_1}:=III_1+III_2.
\label{est-int-1a}
\end{align}
Using  Lemma \ref{sm-ef-int2} and then taking the same way as in the
estimate to \eqref{non-est-1}, we get
\begin{align}
III_1 & \lesssim \sum_{k \in \mathbb{Z}^n, \ |k_1|= k_{\rm max} >4}
\langle k_1 \rangle^{3}  \sum_{\mathbb{S}^{(1)}_{2,1}} \|\Box_k
\left(\Box_{k^{(1)}} u ... \Box_{k^{(\kappa_2+1)}} u \right)
\|_{L^{1}_{x_1} L^2_{x_2,...,x_n}L^2_t(\mathbb{R}^{1+n})}
\nonumber\\
&\lesssim \rho_1(u) (\rho_2(u)+\rho_3(u))^{\kappa_2}.
\label{est-int-1b}
\end{align}
For the estimate of $III_2$, observing the fact that ${\rm supp}
\psi_1 \subset \{\xi: \ |\xi_2| \le 4 |\xi_1|\}$ and using the
multiplier estimate,  then applying \eqref{sm-int-1a}, we have
\begin{align}
III_2 & \lesssim   \sum_{k \in \mathbb{Z}^n, \ |k_1|= k_{\rm max}
>4, \ |k_2| \lesssim |k_1|} \langle k_{1}\rangle^{3}
\sum_{\mathbb{S}^{(1)}_{2,2}} \|\Box_k \left(\Box_{k^{(1)}} u ...
\Box_{k^{(\kappa_2+1)}} u
\right)\|_{L^{(2+\kappa)/(1+\kappa)}_{t,x}(\mathbb{R}^{1+n})
}\nonumber\\
& \lesssim  \rho_3(u)^{1+\kappa_2} . \label{est-int-2}
\end{align}
We need to further control $IV$. Using the decomposition
\eqref{decomp},
\begin{align}
IV & \le \Big \|P_2
\partial^{3}_{x_2} \mathscr{A} \sum_{\mathbb{S}^{(2)}_{2,1}} (\Box_{k^{(1)}} u ...
\Box_{k^{(\kappa_2+1)}} u )  \Big\|_{Y_1} \nonumber\\
& \quad  +  \Big \|P_2
\partial^{3}_{x_2} \mathscr{A} \sum_{\mathbb{S}^{(2)}_{2,2}} (\Box_{k^{(1)}} u ...
\Box_{k^{(\kappa_2+1)}} u )  \Big\|_{Y_1}:=IV_1+IV_2.
\label{est-int-2a}
\end{align}
By Lemma \ref{sm-ef-int2},
\begin{align}
IV_1 & \lesssim \sum_{k \in \mathbb{Z}^n, \ |k_2|= k_{\rm max} >4}
\langle k_2 \rangle^{3}  \sum_{\mathbb{S}^{(2)}_{2,1}} \|\Box_k
\left(\Box_{k^{(1)}} u ... \Box_{k^{(\kappa_2+1)}} u \right)
\|_{L^{1}_{x_1} L^2_{x_2,...,x_n}L^2_t(\mathbb{R}^{1+n})}.
\label{est-int-2b}
\end{align}
By symmetry of $k^{(1)},..., k^{(\kappa_2+1)}$, we can assume that
$|k^{(1)}_2 |=\max_{1\le i \le \kappa_2+1} |k^{(i)}_2 |$ in
$\mathbb{S}^{(2)}_{2,1}$. Using the same way as in the estimate of
$I$, we have
\begin{align}
IV_1 & \lesssim \sum_{\mathbb{S}^{(2)}_{2,1}, \ |k^{(1)}_2| \sim
l_{\rm max} > 4} \langle k^{(1)}_2 \rangle^{3}  \| \Box_{k^{(1)}} u
... \Box_{k^{(\kappa_2+1)}} u
 \|_{L^{1}_{x_1} L^2_{x_2,...,x_n}L^2_t(\mathbb{R}^{1+n})}.
\label{est-int-2ba}
\end{align}
Using H\"older's inequality, we have
\begin{align}
&  \| \Box_{k^{(1)}} u ... \Box_{k^{(\kappa_2+1)}} u
 \|_{L^{1}_{x_1} L^2_{x_2,...,x_n}L^2_t(\mathbb{R}^{1+n})}
 \nonumber\\
& \lesssim  \| \Box_{k^{(1)}} u |\Box_{k^{(2)}} u ...
\Box_{k^{(\kappa_2+1)}} u|^{1/2}
 \|_{L^2_{x,t} (\mathbb{R}^{1+n})} \nonumber\\
  &  \ \ \ \ \times  \||\Box_{k^{(2)}} u ...
\Box_{k^{(\kappa_2+1)}} u|^{1/2}
 \|_{L^{2}_{x_1} L^\infty_{x_2,...,x_n}L^\infty_t(\mathbb{R}^{1+n})} \nonumber\\
& \lesssim  \| \Box_{k^{(1)}} u \|_{L^{\infty}_{x_2}
L^2_{x_1,x_3,...,x_n}L^2_t(\mathbb{R}^{1+n})}
\prod^{\kappa_2+1}_{i=2} \|\Box_{k^{(i)}} u
 \|^{\frac{1}{2}}_{L^{\kappa_2}_{x_2} L^\infty_{x_1,x_3,...,x_n}L^\infty_t(\mathbb{R}^{1+n})
 } \nonumber\\
& \ \ \ \ \times  \prod^{\kappa_2+1}_{i=2} \|\Box_{k^{(i)}} u
 \|^{\frac{1}{2}}_{L^{\kappa_2}_{x_1} L^\infty_{x_2,...,x_n}L^\infty_t(\mathbb{R}^{1+n})
 }.
\label{est-int-2bb}
\end{align}
Observing the inclusion $L^{\kappa}_{x_1}
L^\infty_{x_2,...,x_n}L^\infty_t \bigcap L^{\infty}_{x,t}
 \subset  L^{\kappa_2}_{x_1}
L^\infty_{x_2,...,x_n}L^\infty_t$, we immediately have
\begin{align}
IV_1 & \lesssim  \rho_1(u) (\rho_2(u)+ \rho_3(u))^{\kappa_2}.
\label{est-int-2bc}
\end{align}
For the estimate of $IV_2$,  noticing the fact that ${\rm supp}
\psi_2 \subset \{\xi: \ |\xi_2| \ge 2 |\xi_1|\}$ and applying
\eqref{sm-int-1a}, we have
\begin{align}
IV_2 & \lesssim   \sum_{k \in \mathbb{Z}^n, \ |k_2|= k_{\rm max} >4}
\langle k_2\rangle^{3} \sum_{\mathbb{S}^{(2)}_{2,2}} \|\Box_k
\left(\Box_{k^{(1)}} u ... \Box_{k^{(\kappa_2+1)}} u
\right)\|_{L^{(2+\kappa)/(1+\kappa)}_{t,x}(\mathbb{R}^{1+n})
}\nonumber\\
& \lesssim   \sum_{k \in \mathbb{Z}^n, \ |k_2|= k_{\rm max} >4} \
\sum_{\mathbb{S}^{(2)}_{2,2}} \|\Box_k \left(\Box_{k^{(1)}} u ...
\Box_{k^{(\kappa_2+1)}} u
\right)\|_{L^{(2+\kappa)/(1+\kappa)}_{t,x}(\mathbb{R}^{1+n})
}\nonumber\\
& \lesssim  \rho_3(u)^{1+\kappa_2} . \label{est-int-2c}
\end{align}
The treatment of the other terms in $\rho_1(\cdot)$ is similar.
Therefore, we have shown that
\begin{align}
\rho_1 \left(\mathscr{A} (\sum^n_{i=1} \lambda_i \partial^{3}_{x_i}
u^{\kappa_i+1}) \right) & \lesssim  \sum^n_{i=1} \left(\rho_1(u)
(\rho_2(u)+\rho_3(u))^{\kappa_i} + \rho_3(u)^{1+\kappa_i}\right) .
\label{est-int-2d}
\end{align}

For the  estimate of $\rho_2 (\cdot)$. Denote
\begin{align}
\|u\|_{Z_i} = \sum_{k\in \mathbb{Z}^n} \langle
k\rangle^{3/2-3/\kappa} \|\Box_k u\|_{L^\kappa_{x_i}
L^\infty_{(x_j)_{j\not=i}}L^\infty_t(\mathbb{R}^{1+n})}.
\label{norm-Z}
\end{align}
We have
\begin{align}
\rho_2 \left(\mathscr{A} (\sum^n_{j=1} \lambda_i \partial_{x_i}
u^{\kappa_i+1}) \right) & \lesssim  \sum^n_{i=1}  \left
\|\mathscr{A} (\sum^n_{i=1} \lambda^{3}_i \partial_{x_i}
u^{\kappa_i+1}) \right\|_{Z_j}. \label{est-rho-2a}
\end{align}
Observing the symmetry of $Z_1,...,Z_n$, it suffices to consider the
estimate of $\|\cdot\|_{Z_1}$. Recall that $k_{\rm max}:=\max_{1\le
i \le n} |k_i|$. We have
\begin{align}
\|v\|_{Z_1} & \le \left(\sum_{k\in \mathbb{Z}^n, \, k_{\max}>4} +
\sum_{k\in \mathbb{Z}^n, \, k_{\max} \le 4} \right)\langle
k\rangle^{3/2-3/\kappa} \|\Box_k v\|_{L^\kappa_{x_1}
L^\infty_{x_2,...,x_n}L^\infty_t(\mathbb{R}^{1+n})} \nonumber\\
&:= \Gamma_1(v) + \Gamma_2(v). \label{est-rho-2b}
\end{align}
In view of \eqref{sm-int-2a} and H\"older's inequality,
\begin{align}
 \Gamma_2\left(\mathscr{A} \Big(\sum^n_{i=1} \lambda_i \partial^{3}_{x_i}
u^{\kappa_i+1} \Big) \right)  & \le \sum_{k\in \mathbb{Z}^n, \,
k_{\max} \le 4}
  \left\|\Box_k
\mathscr{A} \Big(\sum^n_{i=1} \lambda_i \partial^{3}_{x_i}
u^{\kappa_i+1} \Big) \right\|_{L^\kappa_{x_1}
L^\infty_{x_2,...,x_n}L^\infty_t(\mathbb{R}^{1+n})} \nonumber\\
& \lesssim  \sum^n_{i=1}
 \sum_{k^{(1)},...,k^{(\kappa_i+1)} \in \mathbb{Z}^n}
\left\|\Box_{k^{(1)}}u ... \Box_{k^{(\kappa_i+ 1)}}u
\right\|_{L^{\frac{2+\kappa_i}{1+\kappa_i}}_{t,x}(\mathbb{R}^{1+n})}
\nonumber\\
& \lesssim  \sum^n_{i=1}
 \sum_{k^{(1)},...,k^{(\kappa_i+1)} \in \mathbb{Z}^n}
\!\!\!\!\!\!\!\!\|\Box_{k^{(1)}}u\|_{L^{2+\kappa_i}_{t,x}(\mathbb{R}^{1+n})}
... \|\Box_{k^{(\kappa_i+ 1)}}u
\|_{L^{2+\kappa_i}_{t,x}(\mathbb{R}^{1+n})} \nonumber\\
& \lesssim  \sum^n_{i=1}  \rho_3(u)^{\kappa_i+1}. \label{est-rho-2c}
\end{align}
It is easy to see that
\begin{align}
\Gamma_1 (v)  & \le \left(\sum_{k\in \mathbb{Z}^n, \, |k_1|=
k_{\max}>4} +...+  \sum_{k\in \mathbb{Z}^n, \, |k_n| = k_{\max} > 4}
\right)\langle k\rangle^{3/2-3/\kappa} \|\Box_k v\|_{L^\kappa_{x_1}
L^\infty_{x_2,...,x_n}L^\infty_t(\mathbb{R}^{1+n})} \nonumber\\
&:= \Gamma^1_1(v) +...+ \Gamma^n_1(v) . \label{est-rho-2d}
\end{align}
Collecting the decomposition \eqref{decomp}, \eqref{sm-int-2a} and
Lemma \ref{sm-ef-int4}, we have
\begin{align}
& \Gamma_1^1 \left(\mathscr{A} \Big(\sum^n_{i=1} \lambda_i
\partial^{3}_{x_i}
u^{\kappa_i+1} \Big) \right) \nonumber\\
& \lesssim \sum^n_{i=1} \sum_{k \in \mathbb{Z}^n, \ |k_1|= k_{\rm
max} >4} \!\!\!\!\!\! \langle k_1 \rangle^{3}
\sum_{\mathbb{S}^{(1)}_{i,1}} \|\Box_k \left(\Box_{k^{(1)}} u ...
\Box_{k^{(\kappa_i+1)}} u \right) \|_{L^{1}_{x_1}
L^2_{x_2,...,x_n}L^2_t(\mathbb{R}^{1+n})}  \nonumber\\
& \ \  + \sum^n_{i=1} \sum_{k \in \mathbb{Z}^n, \ |k_1|= k_{\rm max}
>4} \!\!\! \langle k_1\rangle^{3} \langle k_{\rm
max}\rangle^{3/2}\sum_{\mathbb{S}^{(1)}_{i,2}} \| \Box_k
\left(\Box_{k^{(1)}} u ... \Box_{k^{(\kappa_i+1)}} u
\right)\|_{L^{(2+\kappa_i)/(1+\kappa_i)}_{t,x}(\mathbb{R}^{1+n}) }.
\label{nonhomo-2}
\end{align}
Following the  way as in \eqref{est-I} and \eqref{non-est-2}, we can
get
\begin{align}
\Gamma_1^1 \left(\mathscr{A} \Big(\sum^n_{i=1} \lambda_i
\partial^{3}_{x_i}
u^{\kappa_i+1} \Big) \right) & \lesssim \sum^n_{i=1}\big(\rho_1(u)
(\rho_2(u)+ \rho_3(u))^{\kappa_i} + \rho_3(u)^{1+\kappa_i} \big).
\label{est-lam2}
\end{align}
For the  estimate of $\Gamma^2_1(\cdot)$. By Lemma \ref{sm-ef-int4}
and \eqref{sm-int-2a} ,
\begin{align}
& \Gamma_1^2 \left(\mathscr{A} \Big(\sum^n_{i=1} \lambda_i
\partial^{3}_{x_i}
u^{\kappa_i+1} \Big) \right) \nonumber\\
& \lesssim \sum^n_{i=1} \sum_{k \in \mathbb{Z}^n, \ |k_2|=k_{\max}
>4} \!\!\!\!\!\! \langle k \rangle^{3/2-3/\kappa}
\|\Box_k \left( \mathscr{A}
\partial^{3}_{x_i}
u^{\kappa_i+1}\right) \|_{L^{\kappa}_{x_1}
L^\infty_{x_2,...,x_n}L^\infty_t(\mathbb{R}^{1+n})}  \nonumber\\
& \lesssim \sum^n_{i=1} \sum_{k \in \mathbb{Z}^n, \ |k_2|=k_{\max}
>4} \!\!\!\!\!\! \langle k_2 \rangle^{3}  \sum_{\mathbb{S}^{(2)}_{i,1}}
\|\Box_k \left(\Box_{k^{(1)}} u ... \Box_{k^{(\kappa_i+1)}} u
\right) \|_{L^{1}_{x_2}
L^2_{x_1,x_3,...,x_n}L^2_t(\mathbb{R}^{1+n})}  \nonumber\\
& \ \  + \sum^n_{i=1} \sum_{k \in \mathbb{Z}^n, \ |k_2|=k_{\max}
>4} \!\!\!
\langle k_2\rangle^{3} \langle k_{\rm max}\rangle^{3/2}
\sum_{\mathbb{S}^{(2)}_{i,2}} \| \Box_k \left(\Box_{k^{(1)}} u ...
\Box_{k^{(\kappa_i+1)}} u
\right)\|_{L^{(2+\kappa_i)/(1+\kappa_i)}_{t,x}(\mathbb{R}^{1+n}) } .
\label{nonhomo-2a}
\end{align}
This reduces the same estimate as $\Gamma^1_1(\cdot)$. The terms
$\Gamma^i_1(\cdot)$ for $3\le i \le n$ can be controlled in a
similar way as $\Gamma^2_1(\cdot)$ due to symmetry. Therefore, we
have shown that
\begin{align}
\left\|\mathscr{A} \Big(\sum^n_{i=1} \lambda_i \partial^{3}_{x_i}
u^{\kappa_i+1}\Big) \right\|_{Z_1} & \lesssim  \sum^n_{i=1}
\left(\rho_1(u) (\rho_2(u)+\rho_3(u))^{\kappa_i} +
\rho_3(u)^{1+\kappa_i}\right) . \label{est-Z-8}
\end{align}

Now we turn to estimate $\rho_3(\mathscr{A} \partial^{3}_{x_i}
u^{\kappa_i+1})$. Combining \eqref{st-sm-mo-2} with Lemma
\ref{discret-deriv}, we have
\begin{align}
 \left \|\Box_k \mathscr{A} \partial^{3}_{x_i} f \right\|_{L^\infty_t L^2_x \, \cap \, L^{2+\kappa}_{t,x} (\mathbb{R}^{1+n}) } &
\lesssim    \|\Box_k \partial^{3}_{x_i} f\|_{
L^{(2+\kappa)/(1+\kappa)}_{t,x} (\mathbb{R}^{1+n}) } \nonumber\\
&  \lesssim   \langle k_i\rangle^{3} \|\Box_k  f\|_{
L^{(2+\kappa)/(1+\kappa)}_{t,x} (\mathbb{R}^{1+n}) }.
\label{st-sm-mo-2a}
\end{align}
Collecting \eqref{decomp}, \eqref{st-sm-mo-2a} and
\eqref{st-sm-m-c4}, we have
\begin{align}
& \rho_3(\mathscr{A}
\partial^{3}_{x_1} u^{\kappa_1+1})\nonumber \\
&  \lesssim \sum_{k \in \mathbb{Z}^n, \ |k_1| \le 4} \langle
k_{1}\rangle^{3} \langle k_{\rm max}\rangle^{3/2}
\sum_{k^{(1)},...,k^{(\kappa_1+1)} \in \mathbb{Z}^n } \| \Box_k
\left(\Box_{k^{(1)}} u ... \Box_{k^{(\kappa_1+1)}} u
\right)\|_{L^{(2+\kappa)/(1+\kappa)}_{t,x}(\mathbb{R}^{1+n}) }
\nonumber\\
& + \sum_{k \in \mathbb{Z}^n, \ |k_1| >4} \langle k_1\rangle^{3}
\sum_{\mathbb{S}^{(1)}_{1,1}} \|\Box_k \left(\Box_{k^{(1)}} u ...
\Box_{k^{(\kappa_1+1)}} u \right) \|_{L^{1}_{x_1}
L^2_{x_2,...,x_n}L^2_t(\mathbb{R}^{1+n})} \nonumber \\
& +  \sum_{k \in \mathbb{Z}^n, \ |k_1| >4} \langle k_1\rangle^{3}
\langle k_{\rm max}\rangle^{3/2} \sum_{\mathbb{S}^{(1)}_{1,2}} \|
\Box_k \left(\Box_{k^{(1)}} u ... \Box_{k^{(\kappa_1+1)}} u
\right)\|_{L^{(2+\kappa)/(1+\kappa)}_{t,x}(\mathbb{R}^{1+n}) }.
\label{nonh-1a}
\end{align}
Whether $|k_1|= k_{\rm max} >4 $ or $|k_i|= k_{\rm max} >4,
i=2,...,n$,   using Lemma \ref{sm-ef-int2}, \eqref{est-I} and
\eqref{non-est-2} , we always have
\begin{align}
\rho_3(\mathscr{A} \partial^{3}_{x_1} u^{\kappa_1+1})
 & \lesssim
  \sum^n_{i=1} \left(\rho_1(u) (\rho_2(u)+\rho_3(u))^{\kappa_i} +
\rho_3(u)^{1+\kappa_i}\right). \label{est-rho34}
\end{align}
Until now, we have shown that
\begin{align}
\|\mathscr{T} u\|_X \lesssim \|u_0\|_{M^{3/2}_{2,1}}+ \sum^n_{i=1}
\|u\|^{1+\kappa_i}_X. \label{est-rho34}
\end{align}
Hence, Theorem \ref{DNLS1-mod} holds by a standard contraction
mapping argument. $\hfill \Box$

\section{Proof of Theorem \ref{DNLS-mod}} \label{pf-thm2}

When consider  Theorem \ref{DNLS-mod}, we would like to follow some
ideas as in the proof of Theorem \ref{DNLS1-mod}. However, due to
the nonlinearity contains the general terms $(\partial^{\alpha}_{x}
u)^{\beta}$ with $|\alpha|\leq 3, m+1\leq | \beta | \leq M+1 $, the
proof of Theorem \ref{DNLS1-mod} can not be directly applied.
Inspired by Theorem \ref{DNLS1-mod} , the space $X'$ we need is
likely to be as following:
\begin{align*}
X':= \left\{u\in \mathscr{S}'(\mathbb{R}^{1+n}): \ \|u\|_X:=
\sum^3_{\ell=1} \sum^{3}_{\abs{\alpha}=0} \sum^n_{i=1}
\varrho^{(i)}_\ell(\partial^\alpha_{x} u) \le \delta \right\}.
\end{align*}
where $\alpha$ is a multi-index and
\begin{align*}
& \varrho^{(i)}_1(u) = \sum_{k\in \mathbb{Z}^n, \ |k_i|= k_{\rm
max}>4} \langle k_i\rangle^{3}  \|\Box_k u\|_{L^{\infty}_{x_i}
L^2_{(x_j)_{j\not= i}}L^2_t(\mathbb{R}^{1+n})},\\
& \varrho^{(i)}_2(u) = \sum_{k\in \mathbb{Z}^n} \langle
k\rangle^{3/2-3/m} \|\Box_k u\|_{L^{m}_{x_i}
L^\infty_{(x_j)_{j\not= i}}L^\infty_t(\mathbb{R}^{1+n})},\\
& \varrho^{(i)}_3(u) = \sum_{k\in \mathbb{Z}^n} \langle k
\rangle^{3/2} \|\Box_k u\|_{ L^{2+m}_{x, t} \, \cap \,
L^\infty_tL^2_x (\mathbb{R}^{1+n})}.
\end{align*}

However, when Comparing with the estimates we have established, we
hope the space as following:
\begin{align*}
X:= \left\{u\in \mathscr{S}'(\mathbb{R}^{1+n}): \ \|u\|_X:=
\sum^3_{\ell=1} \sum_{\alpha=0,3} \sum^n_{i,j=1}
\varrho^{(i)}_\ell(\partial^\alpha_{x_j} u) \le \delta \right\}.
\end{align*}
Fortunately, we have  the following:
\begin{lem}\label{norm-eq}
$X$ and $X'$-norm are equivalent.
\end{lem}
{\bf Proof.}
 Obviously, we have $\norm{u}_{X} \leq \norm{u}_{X'}$. For the
 reverse inequality, we only need show that when $\abs{\alpha}=3$ each term in the $X'$ can be controlled by
 $X$.
 Firstly, we consider $\varrho^{(i)}_1(\partial^\alpha_{x} u),
 \abs{\alpha}=3$. Without loss of generality, we can assume $\partial^\alpha_{x} u=\partial_{x_{l}x_{m}x_{o}}u,
  l, m, o=1,...,n$. Then we have
\begin{align*}
 \varrho^{(i)}_{1}(\partial_{x_{l}x_{m}x_{o}}u)= \sum_{k\in \mathbb{Z}^n, \ |k_i|= k_{\rm
max}>4} \langle k_i\rangle^{3}  \|\Box_k
\partial_{x_{l} x_{m} x_{o}}u\|_{L^{\infty}_{x_i} L^2_{(x_j)_{j\not=
i}}L^2_t(\mathbb{R}^{1+n})}
\end{align*}
By symmetry, we can assume $\langle k_{l}\rangle=\max\{\langle
k_{l}\rangle, \langle k_{m}\rangle, \langle k_{o}\rangle\}$ Using
Lemma \ref{discret-deriv} and Remark \ref{eqrm}, we have
\begin{align*}
\|\Box_k \partial_{x_{l}x_{m}x_{o}}u\|_{L^{\infty}_{x_i}
L^2_{(x_j)_{j\not= i}}L^2_t(\mathbb{R}^{1+n})} &\lesssim \langle
k_l\rangle \langle k_m\rangle  \langle k_o\rangle \|\Box_k
u\|_{L^{\infty}_{x_i} L^2_{(x_j)_{j\not=
i}}L^2_t(\mathbb{R}^{1+n})}\\
& \lesssim \langle k_{l}\rangle^{3} \|\Box_k u\|_{L^{\infty}_{x_i}
L^2_{(x_j)_{j\not= i}}L^2_t(\mathbb{R}^{1+n})}\\
& \lesssim \|\Box_k \partial^{3}_{x_{l}}u\|_{L^{\infty}_{x_i}
L^2_{(x_j)_{j\not= i}}L^2_t(\mathbb{R}^{1+n})}
\end{align*}
Hence, we have
\begin{align*}
 \varrho^{(i)}_{1}(\partial_{x_{l}x_{m}x_{o}}u) &\lesssim
 \varrho^{(i)}_{1}(\partial^{3}_{x_{l}}u) \lesssim \|u\|_{X}.
\end{align*}
Secondly, the $\varrho^{(i)}_2(\partial^\alpha_{x} u),
 \abs{\alpha}=3$ can be treat in the same way as above.

Finally, we estimate $\varrho^{(i)}_3(\partial^\alpha_{x} u),
 \abs{\alpha}=3$. Noticing  the denotation of $\varrho^{(i)}_3(\partial^\alpha_{x} u)$, it suffices to show
$\sum_{k\in \mathbb{Z}^n} \langle k \rangle^{3/2} \|\Box_k
u\|_{L^\infty_tL^2_x (\mathbb{R}^{1+n})} \lesssim \norm{u}_{X}$. As
before, we also assume $\partial^\alpha_{x}
u=\partial_{x_{l}x_{m}x_{o}} u$,
  $l, m, o=1,...,n$ and $\langle k_{l}\rangle=\max\{\langle
k_{l}\rangle, \langle k_{m}\rangle, \langle k_{o}\rangle\}$.

Using  Sobolev imbedding Theorem, we have
\begin{align*}
\sum_{k\in \mathbb{Z}^n} \langle k \rangle^{3/2} \|\Box_k
\partial^{\alpha}_{x}u\|_{L^\infty_tL^2_x (\mathbb{R}^{1+n})}& \lesssim \sum_{k\in
\mathbb{Z}^n} \langle k \rangle^{3/2} \|\Box_k
\partial_{x_{l}x_{m}x_{o}}u\|_{L^\infty_tL^2_x (\mathbb{R}^{1+n})}\\
& \lesssim \sum_{k\in \mathbb{Z}^n} \langle k \rangle^{3/2} \|\Box_k
\partial^{3}_{x_{l}}u\|_{L^\infty_tL^2_x (\mathbb{R}^{1+n})}\\
& \lesssim  \varrho^{(i)}_3(\partial^3_{x_l} u)\\ & \lesssim
\norm{u}_{X}.
\end{align*}
 $\hfill\Box$

Considering the following mapping:
\begin{align*}
\mathscr{T}: u(t) \to S(t) u_0 - {\rm i} \mathscr{A}
F((\partial_x^\alpha u)_{\abs{\alpha}\leq 3},
(\partial_x^\alpha\bar{u})_{\abs{\alpha}\leq 3})
\end{align*}
we will show that $\mathscr{T}: X\to X$ is a contraction mapping.

Since $\|u\|_X= \|\bar{u}\|_X$, we can assume that
\begin{align*}
 F((\partial_x^\alpha
u)_{\abs{\alpha}\leq 3},
(\partial_x^\alpha\bar{u})_{\abs{\alpha}\leq 3}) &=
F(\partial_x^\alpha u)_{\abs{\alpha}\leq 3})\nonumber\\
&=\sum_{m+1\le \tilde{R} \le M+1} c_{R_{0},R_{1},R_{2},R_{3}}
u^{R_{0}}(\partial_{x}^{\alpha_{1}}
u)^{R_{1}}(\partial_{x}^{\alpha_{2}}u)^{R_{2}}(\partial_{x}^{\alpha_{3}}u)^{R_{3}},
\end{align*}
where $\tilde{R}=R_{0}+|R_{1}|+|R_{2}|+|R_{3}|$, $R_{i}$ and
$\abs{\alpha_{i}}=i, (i=1,2,3)$ are multi-index.For simiplity, we
denote
\begin{align*}
v_1=...=v_{R_{0}}=u, ... ,
v_{R_{0}+|R_{1}|+|R_{2}|+1}=...=v_{\tilde{R}}=\partial_{x}^{\alpha_{3}}u.
\end{align*}
By \eqref{smo-eff.21}, for $\alpha=0,3$,
\begin{align*}
& \varrho^{(i)}_1(\partial^\alpha_{x_j} S(t)u_0) \lesssim \sum_{k\in
\mathbb{Z}^n, \ |k_i|= k_{\rm max}>4} \langle k_i\rangle^{3/2}
\langle k_j\rangle^{3} \| \Box_k u_0\|_{L^{2} (\mathbb{R}^{n})} \le
\|u_0\|_{ M^{9/2}_{2,1}}.
\end{align*}
By \eqref{st-sm-m-c3}, \eqref{st-sm-m-c1}, we have for $\alpha=0,3$,
\begin{align*}
& \varrho^{(i)}_2(\partial^\alpha_{x_j} S(t)u_0) +
\varrho^{(i)}_3(\partial^\alpha_{x_j} S(t)u_0) \lesssim \|u_0\|_{
M^{9/2}_{2,1}}.
\end{align*}
Hence,
\begin{align*}
\|u\|_X   \lesssim \|u_0\|_{ M^{9/2}_{2,1}}.
\end{align*}
In order to estimate $\varrho^{(i)}_1(\mathscr{A}
\partial^\alpha_{x_j}(v_1... v_{\tilde{R}}))$, $i,j=1,...,n$,  as before it suffices to estimate
 $\varrho^{(1)}_1(\mathscr{A}
\partial^\alpha_{x_1}(v_1... v_{\tilde{R}}))$ and
$\varrho^{(1)}_1(\mathscr{A}
\partial^\alpha_{x_2}(v_1... v_{\tilde{R}}))$.  Similarly as in
\eqref{decomp}, we will use the decomposition
\begin{align}
\Box_k (v_1...v_{\tilde{R}})
 = & \sum_{\mathbb{S}^{(i)}_1} \Box_k
\left(\Box_{k^{(1)}} v_1 ... \Box_{k^{(\tilde{R})}} v_{\tilde{R}}
\right) \nonumber\\
&  +  \sum_{\mathbb{S}^{(i)}_2} \Box_k \left(\Box_{k^{(1)}} v_1 ...
\Box_{k^{(\tilde{R})}} v_{\tilde{R}} \right), \label{decomp2}
\end{align}
where
\begin{align*}
& \mathbb{S}^{(i)}_1:= \{k^{(1)},..., k^{(\tilde{R})} \in \Z^n: \
|k^{(1)}_i|\vee...\vee |k^{(\tilde{R})}_i|>4\}, \\
& \mathbb{S}^{(i)}_2:= \{k^{(1)},..., k^{(\tilde{R})} \in \Z^n: \
|k^{(1)}_i|\vee...\vee |k^{(\tilde{R})}_i|\le 4\}.
\end{align*}

In view of \eqref{1-sm-mod} and \eqref{st-sm-mo-4},
\begin{align}
& \varrho^{(1)}_1(\mathscr{A} \partial^\alpha_{x_1}(v_1...
v_{\tilde{R}})) \nonumber\\
 & \lesssim \sum_{k \in \mathbb{Z}^n, \
|k_1|= k_{\rm max}>4} \langle k_1\rangle^{3}
\sum_{\mathbb{S}^{(1)}_1} \|\Box_k \left(\Box_{k^{(1)}} v_1 ...
\Box_{k^{(\tilde{R})}} v_{\tilde{R}} \right) \|_{L^{1}_{x_1}
L^2_{x_2,...,x_n}L^2_t(\mathbb{R}^{1+n})}  \nonumber\\
& \ \  +  \sum_{k \in \mathbb{Z}^n, \ |k_1|= k_{\rm max}>4} \langle
k_1\rangle^{3} \langle k_{\rm
max}\rangle^{3/2}\sum_{\mathbb{S}^{(1)}_2} \| \Box_k
\left(\Box_{k^{(1)}} v_1 ... \Box_{k^{(\tilde{R})}}
v_{\tilde{R}} \right)\|_{L^{\frac{\tilde{R}+1}{\tilde{R}}}_{t,x}(\mathbb{R}^{1+n}) } \nonumber\\
& := I +II . \label{gnonhomo-1}
\end{align}

Similar to \eqref{est-I},
\begin{align}
I & \lesssim    \sum_{k^{(1)}\in \mathbb{Z}^n, \ |k^{(1)}_1| \sim
k_{\rm max} >4} \langle k^{(1)}_1\rangle  \|\Box_{k^{(1)}}
v_1\|_{L^{\infty}_{x_1} L^2_{x_2,...,x_n}L^2_t(\mathbb{R}^{1+n})} \nonumber\\
& \ \ \times  \sum_{k^{(2)},...,k^{(\tilde{R})}\in \mathbb{Z}^n} \
\prod^{\tilde{R}}_{i=2} \|\Box_{k^{(i)}} v_i
\|_{L^{\tilde{R}-1}_{x_1}
L^\infty_{x_2,...,x_n}L^\infty_t(\mathbb{R}^{1+n})}. \label{gest-I}
\end{align}
By H\"older's inequality and Lemma \ref{qnorm:pnorm},
\begin{align}
&  \|\Box_{k^{(i)}} v_i \|_{L^{\tilde{R}-1}_{x_1}
L^\infty_{x_2,...,x_n}L^\infty_t(\mathbb{R}^{1+n})} \nonumber\\
& \le  \|\Box_{k^{(i)}} v_i \|^{\frac{m}{\tilde{R}-1}}_{L^{m}_{x_1}
L^\infty_{x_2,...,x_n}L^\infty_t(\mathbb{R}^{1+n})} \|\Box_{k^{(i)}}
v_i \|^{1- \frac{m}{\tilde{R}-1}}_{
L^\infty_{x,t}(\mathbb{R}^{1+n}) } \nonumber\\
& \lesssim  \|\Box_{k^{(i)}} v_i
\|^{\frac{m}{\tilde{R}-1}}_{L^{m}_{x_1}
L^\infty_{x_2,...,x_n}L^\infty_t(\mathbb{R}^{1+n})} \|\Box_{k^{(i)}}
v_i \|^{1- \frac{m}{\tilde{R}-1}}_{ L^\infty_t
L^2_{x}(\mathbb{R}^{1+n}) }. \label{gest-I-con}
\end{align}
Owning to Lemma \ref{equivalent} and Remark \ref{eqrm}, it suffices
to consider four cases: $v_i=u$ , $v_i=u_{x_j}$ , $v_i=u_{x_{jj}} $
and $v_{i}=u_{x_{jjj}}, j=1,...,n$. Collecting \eqref{gest-I} and
\eqref{gest-I-con}, we have
\begin{align}
I \lesssim \|u\|_X^{\tilde{R}}. \label{I-control}
\end{align}
Similar to \eqref{non-est-2}, we see that $|k_1| \le C$ in the
summation of $II$. Again, in view of H\"older's inequality and Lemma
\eqref{qnorm:pnorm},
\begin{align}
\|\Box_{k^{(1)}} v_1 ... \Box_{k^{(\tilde{R})}}
v_{\tilde{R}}\|_{L^{\frac{\tilde{R}|+1}{\tilde{R}}}_{x,t}
(\mathbb{R}^{1+n})}  & \le \prod^{\tilde{R}}_{i=1} \|\Box_{k^{(i)}}
v_i
\|_{L^{\tilde{R}+1}_{x,t} (\mathbb{R}^{1+n})} \nonumber\\
& \lesssim \prod^{\tilde{R}|}_{i=1}  \|\Box_{k^{(i)}} v_i
\|_{L^{2+m}_{x,t} \, \bigcap \, L^\infty_t L^2_{x}
(\mathbb{R}^{1+n})}. \label{gest-II-con}
\end{align}

Hence, using a similar way as in \eqref{non-est-2},
\begin{align}
II \lesssim \|u\|_X^{\tilde{R}}. \label{II-control}
\end{align}
We now give the estimate of $\varrho^{(1)}_1(\mathscr{A}
\partial^\alpha_{x_2}(v_1... v_{\tilde{R}}))$.
Since we have obtained the estimate in the case $\alpha=0$, it
suffices to consider the case $\alpha=3$.
 Let $\psi_i$ $(i=1,2)$ be as in Lemma
\ref{sm-ef-int2} and $P_i= \mathscr{F}^{-1} \psi_i \mathscr{F}$. We
have
\begin{align}
& \varrho^{(1)}_1(\mathscr{A}
\partial^{3}_{x_2}(v_1... v_{\tilde{R}})) \nonumber\\
& \le  \sum_{k\in \mathbb{Z}^n, \ |k_1|= k_{\rm max}>4} \langle
k_1\rangle^{3} \|P_1 \Box_k (\mathscr{A}
\partial^{3}_{x_2}(v_1... v_{\tilde{R}}))\|_{L^\infty_{x_1} L^2_{x_2,...,x_n}L^2_t} \nonumber\\
& \ \ + \sum_{k\in \mathbb{Z}^n, \ |k_1|= k_{\rm max}>4} \langle
k_1\rangle^{3} \|P_2 \Box_k (\mathscr{A}
\partial^{3}_{x_2}(v_1... v_{\tilde{R}}))\|_{L^\infty_{x_1} L^2_{x_2,...,x_n}L^2_t} \nonumber\\
& := III+IV. \label{rho112-1}
\end{align}

Using the decomposition \eqref{decomp2},
\begin{align}
III & \le  \sum_{k\in \mathbb{Z}^n, \ |k_1|= k_{\rm max}>4} \langle
k_1\rangle^{3} \sum_{\mathbb{S}^{(1)}_1}   \|P_1 \Box_k (\mathscr{A}
\partial^{3}_{x_2}(\Box_{k^{(1)}} v_1 ... \Box_{k^{(\tilde{R})}}
v_{\tilde{R}} ))\|_{L^\infty_{x_1} L^2_{x_2,...,x_n}L^2_t} \nonumber\\
& \ \ + \sum_{k\in \mathbb{Z}^n, \ |k_1|= k_{\rm max}>4} \langle
k_1\rangle^{3} \sum_{\mathbb{S}^{(1)}_2}   \|P_1 \Box_k (\mathscr{A}
\partial^{3}_{x_2}(\Box_{k^{(1)}} v_1 ... \Box_{k^{(\tilde{R})}}
v_{\tilde{R}}))\|_{L^\infty_{x_1} L^2_{x_2,...,x_n}L^2_t} \nonumber\\
& := III_1 + III_2 \label{rho112-2}.
\end{align}

By Lemma \ref{sm-ef-int2},
\begin{align}
III_1 & \lesssim   \sum_{\mathbb{S}^{(1)}_1} \sum_{k\in
\mathbb{Z}^n, \ |k_1|= k_{\rm max}>4} \langle k_1\rangle^{3}     \|
\Box_k (\Box_{k^{(1)}} v_1 ... \Box_{k^{(\tilde{R})}} v_{\tilde{R}}
)\|_{L^1_{x_1} L^2_{x_2,...,x_n}L^2_t}.
 \label{rho112-3}
\end{align}
By symmetry, we may assume $|k^{(1)}_1|= \max (|k^{(1)}_1|,...,
|k^{(\tilde{R})}_1|)$ in $\mathbb{S}^{(1)}_1$. Hence,
\begin{align}
III_1 & \lesssim   \sum_{\mathbb{S}^{(1)}_1, \ |k^{(1)}_1| \sim
k_{\rm max} > 4} \langle k^{(1)}_1\rangle^{3}   \|\Box_{k^{(1)}}
v_1\|_{L^\infty_{x_1} L^2_{x_2,...,x_n}L^2_t}
\prod^{\tilde{R}}_{i=2} \| \Box_{k^{(i)}} v_i
\|_{L^{\tilde{R}-1}_{x_1}
L^\infty_{x_2,...,x_n}L^\infty_t} \nonumber\\
& \lesssim   \varrho^{(1)}_1 (v_1) \prod^{\tilde{R}}_{i=2}
(\varrho^{(1)}_2 (v_i) + \varrho^{(1)}_3 (v_i)) \lesssim
\|u\|^{\tilde{R}}_X .
 \label{rho112-4}
\end{align}

Applying \eqref{sm-int-1a} and using a similar way as in
\eqref{est-int-2},
\begin{align}
III_2 & \lesssim  \!\!\! \sum_{k \in \mathbb{Z}^n, \ |k_1|= k_{\rm
max}>4, \ |k_2| \lesssim |k_1|} \!\!\! \langle k_1\rangle^{3}
\langle k_{\rm max}\rangle^{3/2} \sum_{\mathbb{S}^{(1)}_2} \|\Box_k
(\Box_{k^{(1)}} v_1 ... \Box_{k^{(\tilde{R})}} v_{\tilde{R}}
)\|_{L^{(2+m)/(1+m)}_{t,x}(\mathbb{R}^{1+n})
}\nonumber\\
& \lesssim  \prod^{\tilde{R}}_{i=1} \varrho^{(1)}_3(v_i) \le
\|u\|^{\tilde{R}}_X. \label{rho112-5}
\end{align}
So, we have shown that
\begin{align}
III & \lesssim  \|u\|^{\tilde{R}}_X. \label{rho112-6}
\end{align}
Now we estimate $IV$. Using the decomposition \eqref{decomp2},
\begin{align}
IV & \le  \sum_{k\in \mathbb{Z}^n, \ |k_1|= k_{\rm max}>4} \langle
k_1\rangle^{3} \sum_{\mathbb{S}^{(2)}_1}   \|P_2 \Box_k (\mathscr{A}
\partial^{3}_{x_2}(\Box_{k^{(1)}} v_1 ... \Box_{k^{(\tilde{R})}}
v_{\tilde{R}} ))\|_{L^\infty_{x_1} L^2_{x_2,...,x_n}L^2_t} \nonumber\\
& \ \ + \sum_{k\in \mathbb{Z}^n, \ |k_1|= k_{\rm max}>4} \langle
k_1\rangle^{3} \sum_{\mathbb{S}^{(2)}_2}   \|P_2 \Box_k (\mathscr{A}
\partial^{3}_{x_2}(\Box_{k^{(1)}} v_1 ... \Box_{k^{(\tilde{R})}}
v_{\tilde{R}}))\|_{L^\infty_{x_1} L^2_{x_2,...,x_n}L^2_t} \nonumber\\
& := IV_1 + IV_2 \label{rho112-7}.
\end{align}
By Lemma \ref{sm-ef-int2},
\begin{align}
IV_1 & \lesssim   \sum_{\mathbb{S}^{(2)}_1} \sum_{k\in \mathbb{Z}^n,
\ |k_2|= k_{\rm max}>4} \langle k_2\rangle ^{3}    \| \Box_k
(\Box_{k^{(1)}} v_1 ... \Box_{k^{(\tilde{R})}} v_{\tilde{R}}
)\|_{L^1_{x_1} L^2_{x_2,...,x_n}L^2_t}.
 \label{rho112-8}
\end{align}
In view of the symmetry, one can bound $IV_1$ by using the same way
as that of $III_1$ and as in
\eqref{est-int-2b}--\eqref{est-int-2bc}:
\begin{align}
IV_1 & \lesssim  \|u\|^{\tilde{R}}_X. \label{rho112-9}
\end{align}
For the estimate of $IV_2$, we apply \eqref{sm-int-1a},
\begin{align}
IV_2 & \lesssim  \!\!\! \sum_{k \in \mathbb{Z}^n, \ |k_1|= k_{\rm
max}>4} \!\!\! \langle k_1\rangle^{3/2} \langle k_2\rangle^{3}
\sum_{\mathbb{S}^{(2)}_2} \|P_2 \Box_k (\Box_{k^{(1)}} v_1 ...
\Box_{k^{(\tilde{R})}} v_{\tilde{R}}
)\|_{L^{(2+m)/(1+m)}_{t,x}(\mathbb{R}^{1+n})
}\nonumber\\
& \lesssim  \sum_{\mathbb{S}^{(2)}_2} \|\Box_{k^{(1)}} v_1 ...
\Box_{k^{(\tilde{R})}} v_{\tilde{R}}
\|_{L^{(2+m)/(1+m)}_{t,x}(\mathbb{R}^{1+n})}\lesssim
\|u\|^{\tilde{R}}_X. \label{rho112-10}
\end{align}
Hence, in view of \eqref{rho112-9} and  \eqref{rho112-10}, we have
\begin{align}
IV \lesssim \|u\|^{\tilde{R}}_X. \label{rho112-11}
\end{align}
Collecting \eqref{I-control}, \eqref{II-control}, \eqref{rho112-6}
and \eqref{rho112-11}, we have shown that
\begin{align}
\sum_{\alpha=0,3} \sum^n_{i,j=1}\varrho_1^{(i)}
(\mathscr{A}\partial^\alpha_{x_j}
(u^{R_{0}}(\partial_{x}^{\alpha_{1}}
u)^{R_{1}}(\partial_{x}^{\alpha_{2}}u)^{R_{2}}(\partial_{x}^{\alpha_{3}}u)^{R_{3}}))
\lesssim \|u\|_X^{\tilde{R}}. \label{lambda1-control}
\end{align}
For later estimate, we need a nonlinear mapping estimate
\begin{lem}[ \cite{Wang2}, Lemma 7.1.]\label{lem5.1}
Let $s\ge 0$,  $1\le p, p_i, \gamma, \gamma_i \le \infty$ satisfy
\begin{align}
\frac{1}{p}= \frac{1}{p_1}+...+\frac{1}{p_N}, \quad
\frac{1}{\gamma}= \frac{1}{\gamma_1}+...+   \frac{1}{\gamma_N}.
\label{p-gamma}
\end{align}
Then
\begin{align}
\sum_{k\in \mathbb{Z}^n} \langle k_1\rangle^{s} \left\|\Box_k (u_1
... u_N) \right\|_{L^{\gamma}_t L^p_x (\mathbb{R}^{1+n})} & \lesssim
\prod^N_{i=1} \left(\sum_{k\in \mathbb{Z}^n} \langle
k_1\rangle^{s}\|\Box_k u_i\|_{L^{\gamma_i}_t L^{p_i}_x
(\mathbb{R}^{1+n})}\right). \label{P7}
\end{align}
\end{lem}
Firstly,  We estimate $\varrho_2^{(1)} (\mathscr{A}
(u^{R_{0}}(\partial_{x}^{\alpha_{1}}
u)^{R_{1}}(\partial_{x}^{\alpha_{2}}u)^{R_{2}}(\partial_{x}^{\alpha_{3}}u)^{R_{3}}))$
and \\
 $\varrho_3^{(1)} (\mathscr{A} (u^{R_{0}}(\partial_{x}^{\alpha_{1}}
u)^{R_{1}}(\partial_{x}^{\alpha_{2}}u)^{R_{2}}(\partial_{x}^{\alpha_{3}}u)^{R_{3}}))$.
 In view of \eqref{st-sm-m-c9} and \eqref{st-sm-mo-6},
\begin{align}
& \sum_{j=2,3} \varrho_j^{(1)} (\mathscr{A}
(u^{R_{0}}(\partial_{x}^{\alpha_{1}}
u)^{R_{1}}(\partial_{x}^{\alpha_{2}}u)^{R_{2}}(\partial_{x}^{\alpha_{3}}u)^{R_{3}})\nonumber\\
& \lesssim  \sum_{k\in \mathbb{Z}^n} \langle k \rangle^{3/2}
\|\Box_k (u^{R_{0}}(\partial_{x}^{\alpha_{1}}
u)^{R_{1}}(\partial_{x}^{\alpha_{2}}u)^{R_{2}}(\partial_{x}^{\alpha_{3}}u)^{R_{3}})\|_{L^{\frac{2+m}{1+m}}_{t,x}(\mathbb{R}^{1+n})
}. \label{lambda23-control}
\end{align}
We use Lemma \ref{lem5.1} to control the right hand side of
\eqref{lambda23-control}:
\begin{align}
& \sum_{k\in \mathbb{Z}^n} \langle k\rangle^{3/2} \|\Box_k
(v_1...v_{\tilde{R}})
\|_{L^{\frac{2+m}{1+m}}_{t,x}(\mathbb{R}^{1+n})} \nonumber\\
& \lesssim \prod^{m+1}_{i=1} \left(\sum_{k\in \mathbb{Z}^n} \langle
k\rangle^{3/2}\|\Box_k v_i\|_{L^{2+m}_{t,x}
(\mathbb{R}^{1+n})}\right)  \prod^{\tilde{R}}_{i=m+2}
\left(\sum_{k\in \mathbb{Z}^n} \langle k \rangle^{1/2}\|\Box_k
v_i\|_{L^{\infty}_{t,x} (\mathbb{R}^{1+n})}\right) \nonumber\\
& \lesssim \prod^{m+1}_{i=1} \left(\sum_{k\in \mathbb{Z}^n} \langle
k \rangle^{3/2}\|\Box_k v_i\|_{L^{2+m}_{t,x}
(\mathbb{R}^{1+n})}\right) \prod^{\tilde{R}}_{i=m+2}
\left(\sum_{k\in \mathbb{Z}^n} \langle k \rangle^{3/2}\|\Box_k
v_i\|_{L^{\infty}_{t}L^2_{x} (\mathbb{R}^{1+n})}\right) \nonumber\\
& \lesssim \prod^{\tilde{R}}_{i=1} \varrho^{(1)}_3(v_i) \le
\|u\|^{\tilde{R}}_X.  \label{lambda23-cont}
\end{align}

Secondly, we estimate $\varrho_2^{(1)} (\mathscr{A}
\partial^{3}_{x_1}(u^{R_{0}}(\partial_{x}^{\alpha_{1}}
u)^{R_{1}}(\partial_{x}^{\alpha_{2}}u)^{R_{2}}(\partial_{x}^{\alpha_{3}}u)^{R_{3}}))$.
\begin{align}
&  \varrho^{(1)}_2 (\mathscr{A} \partial^{3}_{x_1}(v_1...
v_{\tilde{R}})) \nonumber\\
& \lesssim \sum_{k \in \mathbb{Z}^n, \ k_{\max}
>4} \langle
k \rangle^{3/2-3/m}  \|\Box_k \mathscr{A} \partial^{3}_{x_1} \left(
v_1 ... v_{\tilde{R}} \right) \|_{L^{m}_{x_1}
L^\infty_{x_2,...,x_n}L^\infty_t(\mathbb{R}^{1+n})}  \nonumber\\
& \ \ \ \ +  \sum_{k \in \mathbb{Z}^n, \ k_{\max} \le 4} \langle k
\rangle^{3/2-3/m} \|\Box_k \mathscr{A} \partial^{3}_{x_1} \left( v_1
... v_{\tilde{R}} \right) \|_{L^{m}_{x_1}
L^\infty_{x_2,...,x_n}L^\infty_t(\mathbb{R}^{1+n})} \nonumber\\
& := V+VI . \label{gnho-final}
\end{align}

By \eqref{st-sm-mo-6} and Lemma \ref{lem5.1}, we have
\begin{align}
VI \lesssim \sum_{k\in \mathbb{Z}^n}  \|\Box_k (v_1...v_{\tilde{R}})
\|_{L^{\frac{2+m}{1+m}}_{t,x}(\mathbb{R}^{1+n})} \lesssim
\|u\|^{\tilde{R}}_X. \label{gnonhomo-B}
\end{align}
It is easy to see that
\begin{align}
V & \lesssim \left(\sum_{k \in \mathbb{Z}^n, \ |k_1|= k_{\max}
>4}+...+ \sum_{k \in \mathbb{Z}^n, \ |k_n|= k_{\max}>4} \right)
\langle k \rangle^{3/2-3/m} \nonumber\\
& \ \ \ \ \ \times \|\Box_k \mathscr{A}
\partial^{3}_{x_1} \left( v_1 ... v_{\tilde{R}} \right)
\|_{L^{m}_{x_1}
L^\infty_{x_2,...,x_n}L^\infty_t(\mathbb{R}^{1+n})}:=
\Upsilon_1(u)+...+ \Upsilon_n(u). \label{gnho-finala}
\end{align}
Applying the decomposition \eqref{decomp2}, \eqref{sm-int-2a} and
Lemma \ref{sm-ef-int4} , we obtain that
\begin{align}
\Upsilon_1(u) & \lesssim \sum_{k \in \mathbb{Z}^n, \ |k_1|= k_{\rm
max}>4} \langle k_1\rangle^{3}  \sum_{\mathbb{S}^{(1)}_1} \|\Box_k
\left(\Box_{k^{(1)}} v_1 ... \Box_{k^{(\tilde{R})}} v_{\tilde{R}}
\right) \|_{L^{1}_{x_1}
L^2_{x_2,...,x_n}L^2_t(\mathbb{R}^{1+n})}  \nonumber\\
& \ \  +  \sum_{k \in \mathbb{Z}^n, \ |k_1|= k_{\rm max}>4} \langle
k_1\rangle^{9/2} \sum_{\mathbb{S}^{(1)}_2} \| \Box_k
\left(\Box_{k^{(1)}} v_1 ... \Box_{k^{(\tilde{R})}} v_{\tilde{R}}
\right)\|_{L^{\frac{\tilde{R}+1}{\tilde{R}}}_{t,x}(\mathbb{R}^{1+n})
}, \label{gnonhomo-A}
\end{align}
which reduces to the case $\alpha=3$ in \eqref{gnonhomo-1}. So,
\begin{align}
\Upsilon_1(u)  & \lesssim  \|u\|^{\tilde{R}}_X. \label{gnonh-A+Ba}
\end{align}
Again, in view of \eqref{sm-int-2a} and Lemma \ref{sm-ef-int4},
\begin{align}
\Upsilon_2(u) & \lesssim \sum_{k \in \mathbb{Z}^n, \ |k_2|= k_{\rm
max}>4} \langle k_2\rangle^{3}  \sum_{\mathbb{S}^{(2)}_1} \|\Box_k
\left(\Box_{k^{(1)}} v_1 ... \Box_{k^{(\tilde{R})}} v_{\tilde{R}}
\right) \|_{L^{1}_{x_2}
L^2_{x_1,x_3,...,x_n}L^2_t(\mathbb{R}^{1+n})}  \nonumber\\
& \ \  +  \sum_{k \in \mathbb{Z}^n, \ |k_2|= k_{\rm max}>4} \langle
k_2\rangle^{9/2} \sum_{\mathbb{S}^{(2)}_2} \| \Box_k
\left(\Box_{k^{(1)}} v_1 ... \Box_{k^{(\tilde{R})}} v_{\tilde{R}}
\right)\|_{L^{\frac{\tilde{R}+1}{\tilde{R}}}_{t,x}(\mathbb{R}^{1+n})
}, \label{gnonhomo-Aa}
\end{align}
which reduces to the same estimate as $\Upsilon_1(u) $.  Using the
same way as $\Upsilon_2(u)$, we can get the estimates of
$\Upsilon_3(u),...,\Upsilon_n(u)$. So,
\begin{align}
 \varrho^{(1)}_2 (\mathscr{A} \partial^{3}_{x_1}(v_1...
v_{\tilde{R}})) & \lesssim  \|u\|^{\tilde{R}}_X. \label{gnonh-A+B}
\end{align}
We need to further bound $ \varrho^{(1)}_2 (\mathscr{A}
\partial^{3}_{x_i}(v_1... v_{\tilde{R}}))$, $i=2,...,n$, which is essentially the same as $ \varrho^{(1)}_2 (\mathscr{A}
\partial^{3}_{x_1}(v_1... v_{\tilde{R}}))$. Obviously,
  \eqref{gnho-final} holds if we substitute $\partial^{3}_{x_1}$
with $\partial^{3}_{x_i}$. Moreover, using \eqref{sm-int-2a}, Lemmas
\ref{lem5.1} and \ref{sm-ef-int4}, we easily get that
\begin{align}
\rho^{(1)}_2 \left(\mathscr{A} (
\partial^{3}_{x_i} ( v_1... v_{\tilde{R}})) \right) & \lesssim \|u\|_X^{\tilde{R}}. \label{es-rho21-8}
\end{align}
By Lemma \ref{discret-deriv}, \eqref{st-sm-mo-2}, we see that
\begin{align}
\|\Box_k \mathscr{A} \partial^{3}_{x_1} f\|_{L^\infty_t L^2 \, \cap
\,
 L^{2+m}_{t,x}(\mathbb{R}^{1+n})} \lesssim \langle k_1\rangle^{3}
 \|\Box_k f\|_{L^{\frac{2+m}{1+m}}_{t,x}(\mathbb{R}^{1+n})}.
\label{st-sm-m-7a}
\end{align}
Hence, in view of \eqref{st-sm-mo-3}, repeating the procedure as in
the estimates of $\rho_3(u)$ in Theorem \ref{DNLS1-mod},
$\varrho^{(1)}_3 (\mathscr{A}
\partial^{3}_{x_1}(v_1... v_{\tilde{R}}))$ can be controlled by the
right hand side of \eqref{gnonhomo-A} and \eqref{gnonhomo-B}. By
symmetry, the estimate of $\varrho^{(1)}_3 (\mathscr{A}
\partial^{3}_{x_i}(v_1... v_{\tilde{R}})), i=2,...,n$ is identical to $\varrho^{(1)}_3 (\mathscr{A}
\partial^{3}_{x_1}(v_1... v_{\tilde{R}}))$. Summarizing the estimates as in the above, we have shown
that
\begin{align}
\|\mathscr{T} u\|_X \lesssim \|u_0\|_{M^{9/2}_{2,1}} + \sum_{m+1\le
\tilde{R} \le M+1} \|u\|^{\tilde{R}}_X.  \label{est-full}
\end{align}
Therefore, the desired result holds by a standard contracting
mapping argument.

\section*{Acknowledgments}
The author would like to express his great thanks to Professor
Baoxiang Wang and Doctor Zihua Guo for their helpful discussion and
valuable suggestion. The author was partially supported by the
National Natural Science Foundation of China (grants numbers
10571004 and 10621061), the 973 Project Foundation of China, grant
number 2006CB805902.

\begin{thebibliography}{99}
\bibitem{Ben} M. Ben-Artzi, H. Koch and J. C. Saut, Dispersion
estimates for fourth order Schr\"odinger equations, C. R. Acad. Sci.
Paris S\'er. I Math., {\bf 330}(2) (2000) 87--92.

\bibitem{Bergh} J. Bergh and J. L\"ofstr\"om, Interpolation Spaces,
An Introduction, Springer-Verlag, 1976.



\bibitem{Ca}  T. Cazenave, An Introduction to Nonlinear Sch\"odinger
Equations, Textos de M\'etodos Matem\'aticos, 26, Universidade
Federal do Rio de Janeiro, 1996.

\bibitem{Chihara} H. Chihara, Local existence for semilinear
Schr\"odinger equations, Math. Japonica, {\bf 42} (1995) 35--52.



\bibitem{CS} P. Constantin and J. C. Saut, {\rm Local smoothing properties of
dispersive equations}, J. Amer. Math. Soc., {\bf 1} (1988),
413--446.

\bibitem{Dys} K. B. Dysthe, Note on a modification to the
nonlinear Schr\"odinger equation for application to deep water
waves, Proc. R. Soc. Lond., {\bf A369} (1979) 105--114.


\bibitem{Fei2} H. G. Feichtinger, Modulation spaces on locally
compact Abelian group, Technical Report, University of Vienna, 1983.
Published in: ``Proc. Internat. Conf. on Wavelet and Applications",
99--140.  New Delhi Allied Publishers, India, 2003.
{http://www.unive.ac.at/nuhag-php/bibtex/
open\_files/fe03-1\_modspa03.pdf}.


\bibitem{GT} J. Ginibre and Y. Tsutsumi, Uniqueness for the
generalized Korteweg-de Vries equations, SIAM J. Math. Anal., {\bf
20}(1989) 1388--1425.


\bibitem{Guo1} B. L. Guo and B. X. Wang, The global Cauchy problem
 and scattering of solutions for nonlinear Schr\"odinger equations in
 $H^s$. Diff. Integral Eqs., {\bf 15} (9) (2002) 1073--1083.

\bibitem{Hao} C. C. Hao, L. Hsiao and B. X. Wang, Wellposedness
for the fourth order nonlinear Schr\"odinger equations, J. Math.
Anal. Appl., {\bf 320}(1) (2006) 246--265.

\bibitem{Hao1}C. C. Hao, L. Hsiao and B. X. Wang, Wellposedness
of the fourth order nonlinear Schr\"odinger equations in
multi-dimension spaces, J. Math. Anal. Appl., {\bf 328}(1) (2007)
58--83.



\bibitem{Huo} Z. H. Huo and Y. L. Jia, The Cauchy problem for the
fourth-order nonlinear Schr\"odinger equation related to the vortex
filament, J. Diff. Eqns., {\bf 214} (2005) 1--35.

\bibitem{Io-Ke} A. Ionescu and C. E.
Kenig, Low-regularity Schr\"odinger maps, II: Global well posedness
in dimensions $d \ge 3$, Commun. Math. Physics, {\bf 271} (2007),
523--559. arXiv:math/0605209v1.


\bibitem{KPV1} C. E. Kenig, G. Ponce and L. Vega, Oscillatory
integrals and regularity of dispersive equations, Indiana Univ.
Math. J., {\bf 40} (1991) 33--69.

\bibitem{KPV2} C. E. Kenig, G. Ponce and L. Vega, Well-posedness
of the initial value problem for the Korteweg-de Vries equation, J.
Amer. Math. Soc., {\bf 4}(1991) 323--347.

\bibitem{KPV3} C. E. Kenig, G. Ponce and L. Vega, Small solutions to
nonlinear Schr\"odinger equations, Ann. Inst. Henri Poincar\'e,
Nonlinear Anal., {\bf 10} (1993) 255--288.

\bibitem{KPV98} C. E. Kenig, G. Ponce and L. Vega, Smoothing effects
and local existence theory for the generalized nonlinear
Schr\"odinger equations, Invent. Math., {\bf 134} (1998) 489--545.

\bibitem{KPV4} C. E. Kenig, G. Ponce and L. Vega, The Cauchy problem for
quasi-linear Schr\"odinger equations. Invent. Math., {\bf 158}
(2004) 343--388.

\bibitem{KePoVe} C. E.
Kenig, G. Ponce and L. Vega, {\rm Oscillatory integrals and
regularity of dispersive equations,} Indiana Univ. Math. J.,  {\bf
40} \rm (1991), 253--288.

\bibitem{KePoRoVe}  C. E.
Kenig, G. Ponce, C. Rolvent, L. Vega,  The genreal quasilinear
untrahyperbolic Schrodinger equation, Advances in Mathematics {\bf
206} (2006), 402--433.


\bibitem{Seg2003} J. Segata, Well-posedness for the fourth order nonlinear
 Schr\"odinger type equation
 related to the vortex filament, Diff. Integral Eqs., {\bf 16} (7) (2003) 841--864.

\bibitem{Seg2004} J. Segata, Remark on well-posedness for the fourth
order nonlinear Schr\"odinger type equation,  Proc. Amer. Math.
Soc., {\bf 132} (2004) 3559--3568.

\bibitem{SuTo} M. Sugimoto amd N. Tomita, The dilation property of
modulation spaces and their inclusion relation with Besov spaces, J.
Funct. Anal., {\bf 248} (2007), 79--106.

\bibitem{Triebel}H. Tribel, Theory of Function
Spaces, Birkh\"{a}user-Verlag, 1983.




\bibitem{Toft} J. Toft, Continuity properties for modulation spaces,
with applications to pseudo-differential calculus, I. J. Funct.
Anal., {\bf 207} (2004), 399--429.

\bibitem{WangNew}B. X. Wang, L. J. Han and C. Y. Huang, Global
Smooth Effects and Well-Posedness for the  Derivative Nonlinear
Schr\"{o}dinger Equaton with Small Rough Data, arXiv:0808.3098V2.

\bibitem{Wang3}B. X. Wang and C. Y. Huang, Frequency-uniform decomposition method for the generalized BO, KdV and NLS equations, J. Diff.
Eqns., {\bf 239}(2007), 213--250.

\bibitem{Wang2}B. X. Wang and H. Hudzik, The global Cauchy problem for
the NLS and NLKG with small rough data, J. Diff. Eqns., {\bf 232}
(2007) 36-73.


\bibitem{Wang1}B. X. Wang and L. F. Zhao and B. L. Guo, Isometric decomposition
operators, function spaces $E^{\lambda}_{p,q}$ and applications to
nonlinear evolution operators, J. Funct. Anal., {\bf 233} (2006)
1--39.


\bibitem{Zhang} H. Zhang, Global well-posedness and scattering for the fourth order nonlinear
Schr\"{o}dinger equations with small data, arXiv:0809.1512v2.


\end{thebibliography}

\end{document}